\documentclass[12pt]{article}
\usepackage{amsfonts,amssymb,version}

\def\bbb{\mathfrak b}

\def\E{{\mathcal E}}
\def\F{{\mathcal F}}

\def\GG{{\mathbb G}}
\def\ggg{\mathfrak g}

\def\HN{\mathop{\rm HN}\nolimits}

\def\I{{\mathcal I}}
\def\mm{\mathfrak m}

\def\L{{\mathcal L}}
\def\LL{{\mathbb L}} 

\def\OO{{\mathcal O}}

\def\PP{{\mathbb P}}
\def\ppp{\mathfrak p}
\def\Q{\mathbb Q}
\def\QQ{{\mathcal Q}}
\def\pp{\mathfrak p}

\def\Z{\mathbb Z}

\def\deg{\mathop{\rm deg}\nolimits}

\def\id{\mathop{\rm id}\nolimits} 
\def\im{\mathop{\rm im}\nolimits} 
\def\ker{\mathop{\rm ker}\nolimits}

\def\proj{\mathop{\rm Proj}\nolimits}

\def\rank{\mathop{\rm rank}\nolimits}
\def\spec{\mathop{\rm Spec}\nolimits}
\def\Spec{\mathop{\rm Spec}\nolimits}

\def\sym{\mathop{\rm Sym}\nolimits}
\def\toto{\stackrel{\to}{{\scriptstyle \to}}}

\mathchardef\mhyphen="2D

\let\hra\hookrightarrow
\let\ov\overline
\let\un\underline
\let\wh\widehat

\newtheorem{theorem}{Theorem}[section]
\newtheorem{proposition}[theorem]{Proposition}
\newtheorem{lemma}[theorem]{Lemma}

\newtheorem{corollary}[theorem]{Corollary}

\def\stm{\refstepcounter{theorem}\paragraph{\thetheorem}}
\def\rem{\refstepcounter{theorem}\paragraph{Remark \thetheorem}}

\def\proof{\paragraph{Proof}}

\makeatletter \def\l@section{\@dottedtocline{1}{0em}{1.2em}} \makeatother

\begin{document}

\centerline{\large\bf Harder-Narasimhan stacks for principal bundles}

\centerline{\large\bf in higher dimensions and arbitrary characteristics}

\bigskip

\centerline{\bf Sudarshan Gurjar and Nitin Nitsure}

\begin{abstract}
Let $G$ be a split reductive group over a
field $k$ of arbitrary characteristic, chosen suitably. 
Let $X\to S$ 
be a smooth projective morphism of locally noetherian $k$-schemes, 
with geometrically connected fibers. 
We show that for each Harder-Narasimhan type $\tau$ for principal 
$G$-bundles,
all pairs consisting of a principal $G$-bundle on a fiber of $X\to S$ together 
with a given canonical reduction of HN-type $\tau$ form an 
algebraic stack $Bun_{X/S}^{\tau}(G)$ over $S$. 
The forgetful $1$-morphism  
$Bun_{X/S}^{\tau}(G) \to Bun_{X/S}(G)$ to the algebraic stack of 
all principal $G$-bundles on fibers of $X\to S$ is 
a schematic morphism, which is of finite type, separated, 
radicial, and induces an isomorphism on residue fields
of all points of  $Bun_{X/S}^{\tau}(G)$.
It factors via an open substack $Bun_{X/S}^{\ngtr \tau}(G)$ of $Bun_{X/S}(G)$,
inducing a finite morphism $Bun_{X/S}^{\tau}(G) \to Bun_{X/S}^{\ngtr \tau}(G)$.
This is a closed embedding if the Behrend conjecture is satisfied by $G$.

The results of this paper hold in arbitrary characteristic, and in fact it gives better
proofs of the results of [G-N-1] and [G-N-2] which had assumed the characteristic to be zero.
Along the way we give a new proof of the existence of a canonical reduction
over any base field in all dimensions, and we also prove openness of semistability and
semicontinuity of canonical type in a family.
\end{abstract}



\bigskip

\section{Introduction}

Let $G$ be a (connected)
split reductive group over a field $k$ of arbitrary characteristic
together with a split maximal torus $T$ and a Borel $B$ containing it, and let
$\ov{C}$ denote the corresponding closed positive Weyl chamber.
Let $S$ be a locally noetherian $k$-scheme, and let $\pi : X \to S$ be
a smooth projective morphism with geometrically connected fibers together with a
relatively very ample line bundle $\OO_{X/S}(1)$. 
This paper constructs a moduli stack $Bun_{X/S}^{\tau}(G)$ for unstable 
principal $G$-bundles on fibers of $X\to S$ 
which have a given canonical (Harder-Narasimhan) type $\tau\in \ov{C}$.

Our results hold under the following hypothesis $(*)$ about the group $G$.

{\it ${\bf (*)}$ Preservation of semistability under field extensions}:  
If $L/K/k$ are extension fields of $k$, if
$H = P/R_u(P)$ is a Levi quotient group
where $P$ is a parabolic in $G$ that contains $B$ and $R_u(P)$ denotes 
its unipotent radical
and if $E$ a semistable principal $H$-bundle  
on a geometrically irreducible smooth projective curve $C$ over $K$, then 
the base change $E_L$ is a semistable principal $H$-bundle on $C_L$. 

Even though we are interested in principal bundles over
higher dimensional projective varieties, it is enough to make the above hypothesis
only for principal bundles over curves. For as a consequence of the Theorem 12
proved in [Gu], the above hypothesis implies a similar statement where
the curve $C$ is replaced by any geometrically irreducible smooth
projective variety $X$ over $K$ of arbitrary dimension, and the principal bundle $E$ 
is defined over a big open subset of $X$. 

The hypothesis ${\bf (*)}$ is necessary in order to be able to  
define a moduli functor for $G$-bundles of any
given canonical type as functoriality would require the type to be
preserved under pullbacks. 
It is known that ${\bf (*)}$ is a consequence of the Behrend conjecture (see [He]). 
If $G = GL_n$ or $SL_n$,
then the Behrend conjecture is satisfied for all $k$, and   
if $char(k) =0$, then it is satisfied for all reductive $G$.
Moreover, the conjecture always holds for 
classical groups, and holds for exceptional simple groups
whenever $char(k)$ is large enough 
(see Theorem 1 in [He]). 

If $E$ is a principal $G$-bundle on a geometrically irreducible smooth
projective variety $X$ over an extension field $K$ of $k$, then under the hypothesis
$(*)$, there exists a unique `canonical reduction' for $E$,
as first proved by Behrend [Be]
in arbitrary characteristics (also see [Bi-Ho]). The result was already stated in [Ram],  
and is easy to prove in characteristic zero (see [A-B]).
We give a simpler and more conceptual proof of existence in arbitrary characteristics
in Section 4.
The canonical reduction is a reduction of the structure group
of $E$ to a standard parabolic subgroup $P\subset G$,
which satisfies some additional properties.
Such a reduction is the analog for principal bundles of the Harder-Narasimhan
filtration $0 \subset \F_1 \subset \ldots  \subset \F_{\ell} = \F$
for a pure dimensional coherent $\OO_X$-module $\F$ on $X$.
However, for principal
bundles a complicating feature is that a canonical reduction is in general possible
only over a big open subscheme $U\subset X$ that depends on $E$,
where `big' means that $\dim (X-U) \le \dim(X) - 2$.
The existence and uniqueness implies that one can attach a well-defined Harder-Narasimhan
type $\tau$ to each principal bundle $E$ as the type of its canonical reduction,
which is an element of the closed positive Weyl chamber $\ov{C}$. 

Let us now look at the relative case, where we have a family of $G$-bundles
on $X/S$ parameterized by an $S$-scheme $T$, which by definition is
a principal $G$-bundle $E$ on $X_T$. For any $\tau \in \ov{C}$,
we want to define the notion of a family of $G$-bundles of constant type $\tau$.
One can define a family of semistable bundles parameterized by $T$ 
by the requirement that for each $t\in T$, the restriction $E_t = E|_{X_t}$ of $E$
to the schematic fiber over $t$ is semistable. This is the special case $\tau =\tau_{ss}$, which
is the `topological type' of the given family $E_T$ (if $G$ is semisimple, then
necessarily $\tau_{ss} =0$).
Such a point-wise definition works well for $\tau = \tau_{ss}$ only because
semistability of $E_t$ is an open condition on the parameter scheme $T$ of any family.
But in general, such a pointwise definition does not work well as it cannot 
keep track of scheme structures on subsets defined by types of the
underlying space of a parameter scheme, 
and their behaviour under pullbacks.
We now turn to the notion of a relative canonical
reduction for a family.


\pagestyle{myheadings}
\markright{Gurjar and Nitsure: Harder-Narasimhan 
stacks in positive characteristic.}


In the case of pure dimensional coherent
$\OO$-modules, a definition of a relative Harder-Narasimhan filtration was
given in [N-3], and it was proved that it gives a stratification of the
stack of pure coherent sheaves on $X/S$ by locally closed substacks. 
This was the starting point for the series [G-N-1], [G-N-2] and the 
present paper, which respectively consider principal $G$-bundles over curves
in characteristic $0$, principal $G$-bundles over higher dimensional projective
varieties in characteristic $0$, and principal $G$-bundles over
higher dimensional projective varieties in arbitrary characteristic.
In fact, the present paper 
gives better proofs of the main results of [G-N-1] and [G-N-2].

When trying to define a relative reduction for
a family of principal bundles in higher dimensions one immediately faces the
obstacle that even in the absolute case, 
the canonical reduction of a bundle $E$ is defined only over a big open subscheme
$U\subset X$ that depends on $E$. As the techniques of moduli theory
crucially depend on coherence of sheaves and projectivity of schemes,
this represents a problem. In the previous paper [G-N-2] of this series,
we introduced the notion of a relative canonical reduction
$[L,\phi]$ of a given type $\tau$ for a family $E$ of
principal bundles on $X/S$ parameterized by an $S$-scheme $T$. This notion,
recalled in Section 5, works in any characteristic and has built into it a modicum of
$\OO$-coherence that enabled us to use techniques based on projectivity.
For each $\tau \in \ov{C}$, it gives rise to an $S$-groupoid 
$Bun_{X/S}^{\tau}(G)$ which attaches to any $S$-scheme $T$ the groupoid 
whose objects are pairs consisting of a principal $G$-bundle $E$ on
$X_T$ and a relative canonical reduction  
of $E$ of Harder-Narasimhan type $\tau$. Let $Bun_{X/S}(G)$ denote 
the algebraic stack of all $G$-bundles on $X/S$.
The main result of this paper is the following.

\begin{theorem}\label{stack version of main theorem} 
The $S$-groupoid $Bun_{X/S}^{\tau}(G)$ is an 
algebraic stack over $S$. The natural forgetful 
$1$-morphism $Bun_{X/S}^{\tau}(G) \to Bun_{X/S}(G)$
is a schematic morphism of finite type, which is separated, 
universally injective (radicial) and induces an isomorphism on
all residue fields of $Bun_{X/S}^{\tau}(G)$. If $Bun_{X/S}^{\ngtr \tau}(G)$ is the
open substack of $Bun_{X/S}(G)$
where the Harder-Narasimhan type is not strictly greater than $\tau$,
then the $1$-morphism $Bun_{X/S}^{\tau}(G) \to Bun_{X/S}^{\ngtr \tau}(G)$ is
proper hence finite.
Moreover, if the Behrend conjecture holds
for $G$ then $Bun_{X/S}^{\tau}(G) \to Bun_{X/S}^{\ngtr \tau}(G)$
is a closed embedding.
\end{theorem}

Even for $G = GL_n$ or $SL_n$ the above result is new, as the type of semistability
in [N-3] is the Gieseker semistability, while by definition, 
semistability for a principal $GL_n$ or $SL_n$ bundle corresponds to $\mu$-semistability of
the corresponding vector bundle. In the case of $G = GL_n$ or $SL_n$,
the Behrend conjecture holds,
so we can conclude that $Bun_{X/S}^{\tau}(G) \to Bun_{X/S}^{\ngtr \tau}(G)$ is
a closed embedding for such a $G$.

The above theorem can be equivalently re-formulated as follows in purely scheme
theoretic terms, without referring to stacks.

\begin{theorem}\label{family version of main theorem} 
Let $X \to S$ be
a smooth projective morphism with geometrically connected fibers, and 
let $E$ be a principal $G$-bundle on $X$. For any $\tau\in\ov{C}$,
there exists a scheme  $S^{\tau}(E)$ over $S$ 
which has the universal property that for any $S$-scheme $T$, the set 
of all relative canonical reductions of type $\tau$ of the pullback $E_T/X_T/T$
is in a natural bijection with the set of all $S$-morphisms from
$T$ to $S^{\tau}(E)$. Moreover, the
morphism $S^{\tau}(E) \to S$ is of finite type, separated and 
universally injective (radicial), and induces an isomorphism
on all residue fields of points of $S^{\tau}(E)$. If $S^{\ngtr \tau}(E) \subset S$
is the open subscheme where the pointwise Harder-Narasimhan
type is not strictly greater than $\tau$,
then the morphism $S^{\tau}(E) \to S^{\ngtr \tau}(E)$ is proper hence finite.
Moreover, if the Behrend conjecture holds
for $G$ then $S^{\tau}(E) \to S^{\ngtr \tau}(E)$ is a closed embedding.
\end{theorem}

In the course of proving the above theorem, we prove Theorem \ref{semicontinuity of type}
which asserts that the HN-types
are an upper semicontinuous function on the parameter scheme $S$ of any family $E$,
in particular, the locus $S^{\ngtr \tau}(E)$ is an open subscheme of $S$ for any $\tau$.
This uses the openness of semicontinuity (Proposition \ref{openness of semistability}),
which we prove also for rationally defined families (this generality is needed
because canonical reductions are only rationally defined in general). 
The properness of $S^{\tau}(E) \to S^{\ngtr \tau}(E)$ comes from
Theorem \ref{valuative criterion}, which is a valuative criterion for
relative reductions, that may be of more general use.

This paper is organized as follows. The sections 2 and 3 are devoted to setting up
and proving some basic facts about root data and reductive groups for later use. 
These two sections are written in greater detail, keeping in mind those 
young readers who are familiar with the necessary algebraic geometry but 
are not necessarily well-versed in the basics of reductive groups
(an expert reader can begin with section 4, and refer back as needed).
A new proof in all dimensions
of the existence of a canonical reduction over an arbitrary base
field is given in section 4 (see Theorem \ref{existence and uniqueness}).
We prove the uniqueness of a canonical reduction
in dimension $\ge 2$ over an algebraically closed base field,
assuming Behrend's proof of uniqueness over curves. For the convenience of readers, 
we also write down an example (see \ref{Example of big open not everything})
which is well known to experts, which shows
that a canonical reduction may not be possible globally when dimension $\ge 2$,
and so a reduction over a big open subscheme is all that one can have in general.
The section 5 addresses preliminaries on relative canonical reductions,
including a proof of descent and effective descent for these.
Our main results, the  
Theorems \ref{stack version of main theorem}, 
\ref{family version of main theorem}, \ref{valuative criterion} and
\ref{semicontinuity of type} 
are proved in section 6. The proofs make use of the 
relative Proj schemes made from 
Grothendieck's Q-sheaves via the Lemma \ref{projective representation},
which we expect to be of independent interest.

\section{Some useful facts about root data}

The main results in this section are Lemmas \ref{R}, \ref{S},
\ref{characters positive combination of simple roots} and
\ref{chi I alpha as a positive multiple}, which may be well known to 
experts. We have given a cohesive account for use in section 4.

\stm{\bf Definition.} A root datum $(X,\Phi, X^{\vee}, \Phi^{\vee})$ consists of 
finitely generated free abelian groups $X$ and $X^{\vee}$ together with
a bi-additive pairing $\langle \, , \,\rangle : X\times X^{\vee} \to \Z$
which makes $X$ and $X^{\vee}$ the duals of each other, and finite subsets
$\Phi \subset X$, $\Phi^{\vee} \subset X^{\vee}$ together with a bijection
$\phi: \Phi \to \Phi^{\vee}$, under which we denote the image of $\alpha \in \Phi$
by $\alpha^{\vee} \in \Phi^{\vee}$, such that the following conditions are satisfied:
(i) If $\alpha \in \Phi$ 
then $\langle \alpha, \alpha^{\vee}\rangle = 2$,
(ii) for each $\alpha \in \Phi$,
the map $s_{\alpha} : X \to X : x \mapsto x - \langle x, \alpha^{\vee} \rangle \alpha$
maps $\Phi$ to itself, and for 
each $\alpha^{\vee} \in \Phi^{\vee}$,
the map
$s_{\alpha^{\vee}} : X^{\vee} \to X^{\vee} :
y \mapsto y - \langle \alpha, y \rangle \alpha^{\vee}$
maps $\Phi^{\vee}$ to itself.

In the notation, the pairing $\langle \, , \,\rangle : X\times X^{\vee} \to \Z$
and the bijection $\phi: \Phi \to \Phi^{\vee}$ are commonly left out, 
and the root datum is simply denoted by 
$(X,\Phi, X^{\vee}, \Phi^{\vee})$.
The tuple $(X^{\vee}, \Phi^{\vee},X,\Phi)$, 
equipped with the transpose of the pairing and the inverse
of the bijection $\phi$, is again a root datum 
called the dual root datum. The elements of $\Phi$ are called the roots
and those of $\Phi^{\vee}$ are called the co-roots of $(X,\Phi, X^{\vee}, \Phi^{\vee})$.

The subgroup $W$ of $Aut(X)$ generated by all 
maps $s_{\alpha}$, where $\alpha$ varies over $\Phi$, 
is called the Weyl group of the root datum. It 
preserves $\Phi$, and the resulting homomorphism from $W$ to the group
of permutations on $\Phi$ is injective, hence $W$ is finite.
For each root $\alpha$, 
the homomorphisms $s_{\alpha} :X \to X$ and $s_{\alpha^{\vee}}: X^{\vee}\to X^{\vee}$ are
adjoints of each other under the duality pairing of $X$ and $X^{\vee}$, that is,
for all $x\in X$ and $y\in X^{\vee}$ we have $\langle s_{\alpha} (x), y\rangle
= \langle x, s_{\alpha^{\vee}}(y)\rangle$. Equivalently, as $s_{\alpha}$ and
$s_{\alpha^{\vee}}$ are each of order $2$, we have
$\langle s_{\alpha} (x),s_{\alpha^{\vee}}(y)\rangle = \langle x, y\rangle$.
It follows that the dual root datum has the same Weyl group, under
the isomorphism that sends $s_{\alpha} \mapsto s_{\alpha^{\vee}}$
(which is the restriction to $W$ of the isomorphism $Aut(X)\to Aut(X^{\vee}): A \mapsto
{^t}\!A^{-1}$ in matrix terms).

If the only multiples in $\Phi$ 
of each $\alpha \in \Phi$ are $\pm \alpha$, then we say that the root datum is reduced
(unless otherwise stated, we will assume that any root datum which occurs in what
follows is reduced). This condition is self-dual.

\stm\label{Dual decompositions}
{\bf Dual decompositions.} Let $F$ be any field (we will take
$F= \Q$ in what what follows), $V$ a finite dimensional vector space over $F$,
and $V^*$ its dual. Given any subset $C \subset V$, we have a vector
subspace $C^{\perp} \subset V^*$, defined by
$C^{\perp} = \{ f \in V^* \, |\, f(w) =0 \mbox{ for all }w \in C \}$.
If $FC \subset V$ is the vector subspace spanned by $C$ then $C^{\perp} = (FC)^{\perp}
= F(C^{\perp})$.
Given any internal direct sum decomposition $V = \oplus_{s\in S}A_s$
indexed by a set $S$, there exists a unique internal direct sum
decomposition $V^* = \oplus_{s\in S}B_s$, such that under the duality pairing
$\langle \, ,\,\rangle : V\times V^* \to F$ we have
$\langle A_p , B_q\rangle =0$ for any $p,q\in S$ with $p\ne q$.
This means the restricted pairing $\langle \, ,\,\rangle : A_s \times B_s \to F$
is non-degenerate for each $s\in S$
and $B_q = (\oplus_{p\ne q}A_p)^{\perp}$.
We will say that such pairs of decompositions $V = \oplus_{s\in S}A_s$ and
$V^* = \oplus_{s\in S}B_s$ are dual to each other.

\stm\label{Direct sum decomposition-I.}
{\bf Direct sum decomposition-I.}
Let $\Q\Phi \subset \Q\otimes X$ (or
$\Q\Phi^{\vee} \subset \Q\otimes X^{\vee}$)
be the vector subspace over $\Q$ generated by
$\Phi$ (or $\Phi^{\vee}$).
For any root datum $(X, \Phi, X^{\vee}, \Phi^{\vee})$,
we have a dual pair of decompositions
$$\Q\otimes X = \Q\Phi \oplus (\Phi^{\vee})^{\perp} \mbox{ and }
\Q\otimes X^{\vee} = \Q\Phi^{\vee} \oplus \Phi^{\perp}.$$
These decompositions are simultaneously duals of each other 
in two senses: they are dual in the sense of (\ref{Dual decompositions}),
and they get interchanged if the root datum $(X,\Phi, X^{\vee}, \Phi^{\vee})$
is replaced by its dual root datum $(X^{\vee}, \Phi^{\vee}, X, \Phi)$.
A similar comment will apply to the various sets of dual decompositions
associated to root data that will occur in what follows.

\stm\label{global virtual characters, central virtual co-characters}
{\bf Global virtual characters, central virtual co-characters.}
The elements of $X$ (or $\Q \otimes X$)
are called the characters (or virtual characters) and that of $X^{\vee}$
(or $\Q \otimes X^{\vee}$) are called
the co-characters (or virtual co-characters) of the root datum. We call the
elements of $(\Phi^{\vee})^{\perp}$ as global virtual characters, and the
elements of $\Phi^{\perp}$ as central virtual co-characters of the root datum.
Note that $(\Phi^{\vee})^{\perp}$ and $\Phi^{\perp}$ are duals of each other.

\stm
\label{Semi-simplicity.}
{\bf Semisimplicity.} 
A root datum is said to be semisimple if any of the following
equivalent conditions is satisfied. (i) $\Q\Phi = \Q\otimes X$,
(ii) $(\Phi^{\vee})^{\perp} = 0$ (the only global virtual character is $0$),
(iii) $\Q\Phi^{\vee} = \Q\otimes X^{\vee}$,
(iv) $\Phi^{\perp} = 0$ (the only central virtual co-character is $0$).

\stm{\bf Bases, simple roots, simple co-roots.}
There exists a subset
$\Delta \subset \Phi$ (called a base for the root datum) which has the property
that (i) $\Delta$ is linearly independent, and (ii) any root $\alpha \in \Phi$ is
a linear combination $\alpha = \sum_{\beta \in \Delta} n_{\beta}\beta$ with $n_{\beta}\in \Z$,
where either each $ n_{\beta} \ge 0$ or each $ n_{\beta} \le 0$.

The action of the Weyl group $W$ on $\Phi$ takes a base to a base.
The resulting action of $W$ on the set of all bases is transitive and free,
making it a principal $W$-set.

Let a base $\Delta\subset \Phi$ be chosen. The elements of $\Delta$ are called
the simple roots (w.r.t. the chosen $\Delta$).
Their corresponding co-roots $\alpha^{\vee} \in \Phi^{\vee}$
are called the simple co-roots, and they form a subset $\Delta^{\vee} \subset \Phi^{\vee}$
called the co-base corresponding to $\Delta$. It is a base for the dual root datum.

\stm{\bf Fundamental weights $\omega_{\alpha}$, fundamental co-weights
$\omega_{\alpha}^{\vee}$.} Given a base $\Delta$ for a root datum, for
each simple root $\alpha \in \Delta$, there exists a unique element
$\omega_{\alpha} \in \Q\Phi$ with the property that for each $\beta \in \Delta$
we must have 
$$\langle \omega_{\alpha}, \beta^{\vee} \rangle =
\left\{
\begin{array}{ll}
  1 & \mbox{ if } \alpha = \beta, \mbox{ and }  \\
  0 & \mbox{ if } \alpha \ne \beta.
\end{array}
\right.
$$
These elements $\omega_{\alpha} \in \Q\Phi$, as $\alpha$ varies over $\Delta$, are called
the fundamental weights. Dually, the fundamental co-weights 
$\omega_{\alpha}^{\vee} \in \Q\Phi^{\vee}$, as $\alpha$ varies over $\Delta$, are
defined by the property for each $\beta \in \Delta$
we must have 
$$\langle \beta, \omega_{\alpha}^{\vee} \rangle =
\left\{
\begin{array}{ll}
  1 & \mbox{ if } \alpha = \beta, \mbox{ and }  \\
  0 & \mbox{ if } \alpha \ne \beta.
\end{array}
\right.
$$
The fundamental weights $\omega_{\alpha}$ form a $\Q$-linear basis for $Q\Phi$
that is dual to the basis $\Delta^{\vee}$ of $Q\Phi^{\vee}$, and 
the fundamental co-weights $\omega_{\alpha}^{\vee}$ form a $\Q$-linear basis
for $Q\Phi^{\vee}$ that is dual to the basis $\Delta$ of $\Q\Phi$.

\stm\label{Comparison of numbers}
{\bf Comparison of numbers.}
For any $\alpha\in \Delta$ we have
$\langle \alpha, \alpha^{\vee}\rangle =2$.  
For any $\alpha, \beta \in \Delta$ we have $\langle \alpha, \beta^{\vee}\rangle \in \Z$,
and $\langle \alpha, \beta^{\vee}\rangle \le 0$ if $\alpha \ne \beta$.
For any $\alpha, \beta \in \Delta$ we have
$\langle \omega_{\alpha}, \omega_{\beta}^{\vee}\rangle \in \Q_{\ge 0}$, and
$\langle \omega_{\alpha}, \omega_{\alpha}^{\vee}\rangle \in \Q_{> 0}$.
{\bf Caution.} Note that $\langle \alpha, \beta^{\vee}\rangle \ne
\langle \beta, \alpha^{\vee}\rangle$ in general (that is, the Cartan matrix
need not be symmetric).

\stm\label{Weyl chamber}
{\bf The closed positive Weyl chamber $\ov{C} \subset \Q\otimes X^{\vee}$}
consists of all virtual co-characters $x \in \Q\otimes X^{\vee}$
such that $\langle \alpha, x \rangle \ge 0$ for all $\alpha \in \Delta$.
Equivalently, $x\in \ov{C}$ if and only if there exist
$y, z \in \Q\otimes X^{\vee}$ where $y = \sum m_{\alpha}\omega_{\alpha}$ where
$m_{\alpha} \in \Q_{\ge 0}$, and $z \in \Phi^{\perp}$, such that $x = y +z$.
For example, if $G$ is a torus then $\ov{C} = \Q\otimes X^{\vee}$.
The set $\ov{C}$ is a cone, in the sense that it is closed under taking
finite linear combinations of its elements with non-negative rational
coefficients.

\stm\label{Direct sum decomposition-III.}
{\bf Direct sum decomposition-III.} If $I\subset \Delta$ 
then $I \cup \{ \omega_{\beta} \,|\, \beta \in \Delta - I\}$ 
is a linear basis for $\Q\Delta$, and  
$I^{\vee} \cup \{ \omega^{\vee}_{\beta} \,|\, \beta \in \Delta - I\}$ 
is a linear basis for $\Q\Delta^{\vee}$, where
$I^{\vee} \subset \Delta^{\vee}$ consists of all $\alpha^{\vee}$ where
$\alpha \in I$. Also note that $\Q\Delta = \Q\Phi \subset
\Q\otimes X$ and $\Q\Delta^{\vee} = \Q\Phi^{\vee} \subset
\Q\otimes X^{\vee}$.
Consequently, we have internal direct sum decompositions
\begin{eqnarray*}
\Q\Delta & = & \Q I \oplus (\oplus_{\beta \in \Delta -I} \Q\omega_{\beta}),\\
(I^{\vee})^{\perp} & = & (\oplus_{\beta \in \Delta -I}\Q\omega_{\beta})
                           \oplus (\Delta^{\vee})^{\perp}, \\
\Q \otimes X & = & \Q\Delta \oplus (\Delta^{\vee})^{\perp} \\
             & = &      \Q I \oplus (I^{\vee})^{\perp} \\
             & = & \Q I \oplus (\oplus_{\beta \in \Delta -I}
                  \Q\omega_{\beta})\oplus (\Delta^{\vee})^{\perp}.\\
\mbox{ Dually},~~~~ & & \\
\Q\Delta^{\vee} & = & \Q I^{\vee} \oplus
                     (\oplus_{\beta \in \Delta -I}\Q\omega^{\vee}_{\beta}), \\
  I^{\perp}       & = & (\oplus_{\beta \in \Delta -I} \Q\omega^{\vee}_{\beta})
                        \oplus \Delta^{\perp},\\ 
\Q\otimes X^{\vee} & = & \Q\Delta^{\vee} \oplus \Delta^{\perp} \\
                  & = & \Q I^{\vee} \oplus I^{\perp} \\
                  & = &   \Q I^{\vee} \oplus
                        (\oplus_{\beta \in \Delta -I} \Q\omega^{\vee}_{\beta})
                        \oplus \Delta^{\perp}. 
\end{eqnarray*}

\stm\label{in roots terms character on parabolic trivial on center of G}
{\bf Characterization of $\oplus_{\beta \in \Delta -I}\Q\omega_{\beta}\subset \Q \otimes X$.}
As an immediate consequence of the decomposition 
(\ref{Direct sum decomposition-III.}), we have
$$\oplus_{\beta \in \Delta -I}\Q\omega_{\beta} =
\Q \Delta \cap (I^{\vee})^{\perp} \subset \Q \otimes X.$$
In other words, suppose that $\chi \in \Q \otimes X$.
Then $\chi \in \oplus_{\beta \in \Delta -I}\Q\omega_{\beta}$ if and only if
$\chi \in \Q \Delta$ and $\langle \chi,  I^{\vee} \rangle =0$. 
Dually, $\oplus_{\beta \in \Delta -I}\Q\omega^{\vee}_{\beta} =
\Q \Delta^{\vee} \cap I^{\perp} \subset \Q \otimes X^{\vee}$.

The following statement indirectly occurs in [Bi-Ho]. 

\begin{lemma}\label{R}
Each $\alpha \in \Delta - I$ can be uniquely expressed as a linear combination 
$$\alpha =   \sum_{\beta \in I} b_{\beta}\beta + \sum_{\gamma \in \Delta -I}
c_{\gamma}\omega_{\gamma},$$
where $b_{\beta}, c_{\gamma} \in \Q$.
Moreover, we must have $c_{\alpha} \ne 0$, and each $b_{\beta} \le 0$.
\end{lemma}

\proof The decomposition
(\ref{Direct sum decomposition-III.}) shows that $\alpha$ can be uniquely
expressed as a linear combination of the roots $\beta \in I$ and the weights
$\omega_{\gamma}$ for $\gamma \in \Delta - I$, so we can uniquely express
$\alpha = \sum_{\beta \in I} b_{\beta}\beta + \sum_{\gamma \in \Delta -I}
c_{\gamma}\omega_{\gamma}$. The set of vectors
$I \cup \{ \alpha \} \cup \{ \omega_{\gamma}\,|\, \gamma \in \Delta -I  - \{ \alpha \} \}$
is linearly independent,
which shows that $c_{\alpha} \ne 0$. Next we prove that 
each $b_{\beta} \le 0$. Note that we can uniquely express each $b_{\beta}$ as
$b'_{\beta} - b''_{\beta}$ where $b'_{\beta}\ge 0$, $b''_{\beta} \ge 0$ and 
$b'_{\beta}b''_{\beta} =0$. Let $A = \sum_{\beta \in I} b_{\beta}\beta$,
$A' = \sum_{\beta \in I} b'_{\beta}\beta$, and $A'' = \sum_{\beta \in I} b''_{\beta}\beta$.
Let $\Omega = \sum_{\gamma \in \Delta -I}
c_{\gamma}\omega_{\gamma}$. Then we have $A = A' - A''$, and hence
$$\alpha = A' - A'' + \Omega.$$
Next, the Weyl group $W$ is a finite subgroup of $Aut(X)\subset GL(\Q\otimes X)$,
and it preserves the decomposition $\Q\otimes X = \Q\Phi \oplus (\Phi^{\vee})^{\perp}$.
Hence we can choose a positive definite inner product
$$(\,,\,) : \Q\otimes X \times \Q\otimes X \to \Q$$
that is preserved by $W$, and moreover, 
under which $\Q\Phi$ is orthogonal to  $(\Phi^{\vee})^{\perp}$.

It is a basic property of such an inner product (called a Killing form)
that for any two distinct simple
roots $\beta \ne \gamma$ in $\Delta$
we have
$$(\beta, \gamma) \le 0, ~~(\beta, \omega_{\gamma}) = 0 ~\mbox{ and }
(\omega_{\beta}, \omega_{\gamma}) \ge 0.$$
Moreover, $(\beta, \omega_{\beta} ) >0$ for each $\beta \in \Delta$.

It follows that $(\alpha, A') \le 0$, $(A'',A') \le 0$ and $(\Omega, A') = 0$.
Hence we get
$$0 \ge (\alpha, A') = (A',A') - (A'',A') + (\Omega, A')
=  (A',A') - (A'',A') \ge 0.$$
As $(A',A') \ge 0$, the above shows that $(A',A') =0$, and hence $A'=0$
by the positive definiteness of the Killing form.
\hfill $\square$

Raghunathan has informed us that the following statement is also
a consequence of Lemma 1.1 in [Rag].

\begin{lemma}\label{S}
For each $\alpha \in I$, the fundamental weight
$\omega_{\alpha}$ can be uniquely expressed as a linear combination 
$$\omega_{\alpha} =   \sum_{\beta \in I} b_{\beta}\beta + \sum_{\gamma \in \Delta -I}
c_{\gamma} \omega_{\gamma},$$
where $b_{\beta}, c_{\gamma} \in \Q$.
Moreover, we must have
$b_{\alpha} >0$, each $b_{\beta} \ge 0$, and each $c_{\gamma} \ge 0$. 
\end{lemma}

\proof The decomposition
(\ref{Direct sum decomposition-III.}) shows that $\omega_{\alpha}$ can be uniquely
expressed as a linear combination
$\sum_{\beta \in I} b_{\beta}\beta + \sum_{\gamma \in \Delta -I}
c_{\gamma} \omega_{\gamma}$
of the roots $\beta \in I$ and the weights
$\omega_{\gamma}$ for $\gamma \in \Delta - I$.
By linear independence of the set of vectors 
$$(I - \{ \alpha \})\cup \{ \omega_{\alpha} \}
\cup \{ \omega_{\gamma}\,|\, \gamma \in \Delta -I\},$$
it follows that $b_{\alpha} \ne 0$.

As in the proof of Lemma \ref{R}, we write each $b_{\beta}$ as
a difference $b'_{\beta}- b''_{\beta}$, and write $A = A' - A''$.
Similarly, we can write each $c_{\gamma}$ uniquely as
a difference $c'_{\gamma}-c''_{\gamma}$ where
$c'_{\gamma}\ge 0$, $c''_{\gamma}\ge 0$ and $c'_{\gamma}c''_{\gamma} =0$.
Let $\Omega' = \sum_{\gamma \in \Delta -I}
c'_{\gamma}\omega_{\gamma}$ and let 
$\Omega'' = \sum_{\gamma \in \Delta -I}
c''_{\gamma}\omega_{\gamma}$, so that we get $\Omega = \Omega ' - \Omega ''$.
Hence we get
$$\omega_{\alpha} = A' - A'' + \Omega' - \Omega''$$
Now applying $(-,A'')$ to both sides of the above equation, we get
$$0 \le (\omega_{\alpha}, A'') = (A',A'') - (A'',A'') \le 0$$
and hence $A'' = 0$. This finishes the proof that each $b_{\beta} \ge 0$.

Hence we get $\omega_{\alpha} = A' + \Omega' - \Omega''$.
The proof will be complete if we show that $\Omega'' = 0$.
Suppose $\Omega'' \ne 0$. Then we can choose $\gamma \in \Delta - I$
such that $c'_{\gamma} = 0$ and $c''_{\gamma} >0$. Applying $(-,\gamma)$
to both sides of the equation $\omega_{\alpha} = A' + \Omega' - \Omega''$
gives
$$0 \le (\omega_{\alpha}, \gamma)
= (A', \gamma) +(\Omega', \gamma) - (\Omega'', \gamma)$$
As $(A', \gamma) \le 0$, $(\Omega', \gamma) = 0$ and $(\Omega'', \gamma) >
c''_{\gamma} (\omega_{\gamma} ,\gamma) > 0$, this gives a contradiction.
Hence $\Omega''=0$, completing the proof.
\hfill $\square$

\stm\label{definition of partial order}
{\bf Partial order on $\Q \otimes X$ and on $\Q \otimes X^{\vee}$.}
For $x\in \Q \otimes X$, we say that $0 \le x$ if we have a rational
linear combination  $x = \sum_{\alpha \in \Delta} m_{\alpha}\alpha$ where each
$m_{\alpha} \ge 0$. We will say that an $x\ge 0$ (or $x>0$)
is a non-negative (or positive) virtual character, 
For $x, y \in \Q \otimes X$, we say that $x\le y$ if
$0\le y - x$. Note that $x \ge 0$ if and only if
$x \in \Q\Delta$ and 
$\langle x, \omega^{\vee}_{\alpha} \rangle \ge 0$ for 
all $\alpha \in \Delta$. For example, any simple root, or any fundamental weight
is a positive character. 
For $u \in \Q \otimes X^{\vee}$, we say that $0 \le u$ if we have a rational
linear combination  
$u = \sum_{\alpha \in \Delta} m_{\alpha}\alpha^{\vee}$ where each
$m_{\alpha} \ge 0$. We will say that an element 
$u\ge 0$ (or $u>0$)
is a non-negative (or positive) virtual co-character. 
For $u, v \in \Q \otimes X_*(T)$, we say that $u\le v$ if
$0\le v - u$. Note that $u \ge 0$ if and only if $u\in \Q \Delta^{\vee}$ and 
$\langle \omega_{\alpha}, u \rangle \ge 0$ for 
all $\alpha \in \Delta$. For example, any simple co-root, or any fundamental 
co-weight is a positive co-character.

\stm\label{basic positive virtual characters for I}
{\bf Basic positive virtual characters associated to any $I\subset \Delta$.}
By statement (\ref{Direct sum decomposition-III.}),
each $\alpha \in \Delta -I$ can be uniquely expressed as
$\alpha =   \sum_{\beta \in I} b_{\beta}\beta + \sum_{\gamma \in \Delta -I}
c_{\gamma} \omega_{\gamma}$
where $b_{\beta}, c_{\gamma} \in \Q$. By Lemma \ref{R}, each $b_{\beta} \le 0$ and
$c_{\alpha}\ne 0$.
Let
$$\chi_{I, \alpha} = \alpha - \sum_{\beta \in I} b_{\beta}\beta
= \sum_{\gamma \in \Delta -I}
c_{\gamma} \omega_{\gamma}.$$
We will call these elements $\chi_{I, \alpha}  \in \Q\otimes X$,
where $\alpha$ varies over $\Delta -I$, as
the basic positive virtual characters associated to $I\subset \Delta$.
Note that $\chi_{I, \alpha}>  0$ in terms of
the partial order on $\Q\otimes X$ defined in statement \ref{definition of partial order}.


\begin{lemma}\label{characters positive combination of simple roots}
For any $I\subset \Delta$, the basic positive virtual characters 
$\chi_{I, \alpha}$ (where $\alpha$ varies over $\Delta -I$)
defined above have the following properties.

(1) The elements $\chi_{I, \alpha}$, as $\alpha$ varies over $\Delta - I$, form a
$\Q$-linear basis
for $\oplus_{\beta \in \Delta - I}\Q\omega_{\beta}$.

(2) An element $\chi \in \oplus_{\beta \in \Delta - I}\Q\omega_{\beta}$
satisfies $\chi \ge 0$ w.r.t. the partial order
\ref{definition of partial order} on $\Q \otimes X$ 
if and only if $\chi$ is a $\Q$-linear combination
$\chi = \sum_{\alpha \in \Delta - I}n_{\alpha}\chi_{I, \alpha}$ 
with each $n_{\alpha}\ge 0$.
\end{lemma}

\proof As
the simple roots $\alpha \in \Delta -I$ are linearly independent,
it follows that the vectors $\chi_{I, \alpha} = \alpha - \sum_{\beta \in I} b_{\beta}\beta$
are linearly independent. These vectors $\chi_{I, \alpha}$ lie in
$\oplus_{\beta \in \Delta - I}\Q\omega_{\beta}$ which is of dimension
equal to the cardinality of $\Delta - I$. Hence the $\chi_{I, \alpha}$ form
a linear basis for the $\Q$-vector space 
$\oplus_{\beta \in \Delta - I}\Q\omega_{\beta}$ as claimed.

We next prove the assertions about positivity.
For each $\alpha \in \Delta - I$, the vector
$\chi_{I, \alpha} = \alpha - \sum_{\beta \in I} b_{\beta}\beta$
is a non-negative linear combination of simple roots, where
by a `non-negative linear combination' of some vectors in a $\Q$-vector space,
we will mean a linear combination all whose coefficients are non-negative
rational numbers.
Hence any non-negative linear combination of the $\chi_{I, \alpha}$'s (where
$\alpha$ varies over $\Delta - I$) is again a non-negative linear
combination of simple roots.

Conversely, suppose an element of $\chi \in \oplus_{\beta \in \Delta - I}\Q\omega_{\beta}$
can be expressed as a linear combination 
$\chi = \sum_{\gamma \in \Delta} p_{\gamma}\gamma$, with each $p_{\gamma} \ge 0$. Let
$$\chi ' = \sum_{\alpha \in \Delta - I}p_{\alpha}\chi_{I, \alpha}.$$
Then note that
$$\chi - \chi' \in \oplus_{\gamma \in I}\Q\gamma.$$
But $\chi, \chi_{I, \alpha} \in \oplus_{\beta \in \Delta - I}\Q\omega_{\beta}$,
and by the direct sum decomposition (\ref{Direct sum decomposition-III.}), we have
$$(\oplus_{\gamma \in I}\Q\gamma)\cap
(\oplus_{\beta \in \Delta - I}\Q\omega_{\beta}) = 0.$$
Hence $\chi = \chi'$ which is a non-negative linear combination of the $\chi_{I, \alpha}$
as desired.
\hfill $\square$

\stm\label{root datum associated to a subset}
{\bf The root datum associated to a subset $I$ of $\Delta$.}
Let $(X, \Phi, X^{\vee}, \Phi^{\vee})$ be a root datum, and let $\Delta\subset \Phi$
be a base. 
To any $I\subset \Delta$, we now associate a root datum
$$(X, \Phi_I, X^{\vee}, \Phi^{\vee}_I)$$
as follows.
The groups $X$, $X^{\vee}$ and the pairing
$\langle \, ,\,\rangle: X\times X^{\vee} \to \Z$ is as before. 
The subset $\Phi_I \subset \Phi$ consists of all $\alpha \in \Phi$
that can be expressed as an integral linear combination of elements of $I$.
Let 
$I^{\vee}$ consist of all $\alpha^{\vee}$ corresponding to $\alpha \in I$,
and let $\Phi^{\vee}_I \subset \Phi^{\vee}$ consists of all 
$\alpha^{\vee} \in \Phi^{\vee}$
that can be expressed as an integral linear combination of elements of $I^{\vee}$.
The bijection $\Phi \to \Phi^{\vee}$ restricts to a bijection
$\Phi_I \to \Phi^{\vee}_I$. The resulting tuple
$(X, \Phi_I, X^{\vee}, \Phi^{\vee}_I)$ can be seen to be a
root datum.

The Weyl group $W_I$ for this root datum is the subgroup of the original Weyl
group $W$ generated by all $s_{\alpha}$ for $\alpha \in \Phi_I$. 

The subset $I\subset \Phi_I$ is a base for this root datum, with $I^{\vee}$
the corresponding co-base. In terms of the decomposition 
$\Q\otimes X = \Q I \oplus (\oplus_{\beta \in \Delta -I}
\Q\omega_{\beta})\oplus (\Delta^{\vee})^{\perp}$ of statement 
(\ref{Direct sum decomposition-III.}), we have
$$(I^{\vee})^{\perp} = (\oplus_{\beta \in \Delta -I}
\Q\omega_{\beta})\oplus (\Delta^{\vee})^{\perp}.$$
Dually, in terms of the decomposition $\Q \otimes X^{\vee} = 
\Q I^{\vee} \oplus (\oplus_{\beta \in \Delta -I}\Q\omega^{\vee}_{\beta})
\oplus \Delta^{\perp}$, we have
$I^{\perp} = (\oplus_{\beta \in \Delta -I}
\Q\omega^{\vee}_{\beta})\oplus \Delta^{\perp}.$
We will denote the fundamental weights and co-weights
for $(X, \Phi_I, X^{\vee}, \Phi^{\vee}_I)$ by
$\omega_{I,\alpha} \in \Q\otimes X$ and $\omega^{\vee}_{I,\alpha} \in \Q\otimes X^{\vee}$,
where $\alpha \in I$. \\
{\bf Caution.} In general,
the elements $\omega_{I,\alpha}$ and $\omega^{\vee}_{I,\alpha}$ (where $\alpha \in I$)
are quite distinct from the fundamental weights
$\omega_{\alpha}$ and $\omega^{\vee}_{\alpha}$. 

\stm\label{positivity of fundamental wts for I w.r.t. Delta}
{\bf Positivity of $\omega_{I,\alpha}$ w.r.t. $\le_{\Delta}$.} Note that the partial order
$\le_{\Delta}$
on $\Q\otimes X$ defined by the base $\Delta$ for the
root datum $(X, \Phi, X^{\vee}, \Phi^{\vee})$
(see \ref{definition of partial order})
differs from the partial order $\le_I$ 
on $\Q\otimes X$ defined by the base $I$ for the
root datum $(X, \Phi_I, X^{\vee}, \Phi^{\vee}_I)$ if $I \ne \Delta$.
But we have an implication
$$ x \le_I y \Rightarrow x \le_{\Delta} y$$
and hence the elements $\omega_{I,\alpha} >_{\Delta} 0$ for $\alpha \in I$.
Dually, $\omega_{I,\alpha}^{\vee} >_{\Delta} 0$ in $\Q \otimes X^{\vee}$ for $\alpha \in I$.

\stm\label{chi I alpha is a global virtual character for the root datum for I} 
{\bf $\chi_{I,\alpha}$ is a global virtual character for the root datum
$(X, \Phi_I, X^{\vee}, \Phi^{\vee}_I)$.}
Note that for the root datum $(X, \Phi_I, X^{\vee}, \Phi^{\vee}_I)$,
the vector spaces of global virtual characters and
central virtual co-characters 
(see statement \ref{global virtual characters, central virtual co-characters})
are $(\Phi^{\vee}_I)^{\perp} = (I^{\vee})^{\perp}$,
and $(\Phi_I)^{\perp} = I^{\perp}$, respectively. Hence, as
for any $\alpha \in \Delta - I$ the basic positive character
virtual character $\chi_{I,\alpha}$ lies in
$\oplus_{\beta \in \Delta - I}\omega_{\beta} \subset (I^{\vee})^{\perp}$,
it follows that $\chi_{I,\alpha}$ is a global virtual character for
the root datum $(X, \Phi_I, X^{\vee}, \Phi^{\vee}_I)$.

\begin{lemma}
\label{chi I alpha as a positive multiple}
{\bf $\chi_{I, \alpha}$ is a positive multiple of $\omega_{I\cup \{ \alpha \}, \alpha}$.}
Let $I \subset \Delta$ and let $\alpha \in \Delta - I$. Let
$\chi_{I, \alpha}$ be the basic positive virtual character
as defined in statement \ref{basic positive virtual characters for I}. 
Let $\omega_{I\cup \{ \alpha \}, \alpha}$ be the fundamental weight
corresponding to the simple root $\alpha \in I\cup \{ \alpha \}$ of the
root datum $(X, \Phi_{I\cup \{ \alpha \}}, X^{\vee}, \Phi^{\vee}_{I\cup \{ \alpha \}})$
with base $I\cup \{ \alpha \}$.
Then  $\langle \chi_{I,\alpha} , \alpha^{\vee} \rangle >0$
and
$$\chi_{I,\alpha} =  \langle \chi_{I,\alpha} , \alpha^{\vee}\rangle 
\omega_{I\cup \{ \alpha \}, \alpha}.$$
\end{lemma}

\proof By it definition,
$\chi_{I, \alpha} \in \Q(I \cup \{ \alpha \} ) \cap
(\oplus_{\gamma \in \Delta - I}\Q\omega_{\gamma})$.
For $\gamma \in \Delta - I$, note that $\omega_{\gamma} \in (I^{\vee})^{\perp}$.
Hence
$$~~~~~~~~~~~
\chi_{I, \alpha} \in \Q(I \cup \{ \alpha \} ) \cap (I^{\vee})^{\perp}
~~~~~~~~~~~~~~~~\ldots (*)$$
As given by statement \ref{root datum associated to a subset}, we have a root datum
$(X, \Phi_{ I \cup \{ \alpha \} }, X^{\vee}, \Phi^{\vee}_{ I \cup \{ \alpha \} })$
associated to the subset $I\cup \{ \alpha \} \subset \Delta$, 
with base $I\cup \{ \alpha \}$.
Let 
$$\omega_{I\cup \{ \alpha \}, \alpha} \in \Q \Phi_{ I \cup \{ \alpha \} }$$ 
denote the fundamental weight associated to its simple root 
$\alpha \in I \cup \{ \alpha \}$.
It follows by statement \ref{character on parabolic trivial on center of G}
applied to the root datum
$(X, \Phi_{ I \cup \{ \alpha \} }, X^{\vee}, \Phi^{\vee}_{ I \cup \{ \alpha \} })$
and the subset $I\subset I\cup \{ \alpha \}$
that
$$~~~~~~~~~~~
\Q \omega_{I\cup \{ \alpha \}, \alpha} = \Q (I\cup \{ \alpha \}) \cap (I^{\vee})^{\perp}
~~~~~~~~~~~~~~~~\ldots (**)$$
From $(*)$ and $(**)$, it follows that 
$\chi_{I, \alpha} \in \Q \omega_{I\cup \{ \alpha \}, \alpha}$. 
As both $\chi_{I,\alpha}$ and $\omega_{I\cup \{ \alpha \}, \alpha}$ are non-negative linear 
combinations of $I\cup  \{ \alpha \}$, and as both are non-zero,
it now follows that
there exists a rational number $a_{\alpha} >0$ such that 
$\chi_{I,\alpha} = a_{\alpha} \omega_{I\cup \{ \alpha \}, \alpha}$.
Applying $\alpha^{\vee}$ to both sides of the above shows that
$a_{\alpha} = \langle \chi_{I,\alpha} , \alpha^{\vee} \rangle$.
Hence $\langle \chi_{I,\alpha} , \alpha^{\vee} \rangle >0$
and we get the formula
$\chi_{I,\alpha} =  \langle \chi_{I,\alpha} , \alpha^{\vee}\rangle 
\omega_{I\cup \{ \alpha \}, \alpha}$.
\hfill $\square$

\stm{\bf Caution.} In general, $\chi_{I,\alpha}$ is {\it not} a scalar multiple of
the fundamental weight $\omega_{\alpha}$ of $(X,\Phi, X^{\vee}, \Phi^{\vee})$.
For example, if $I = \emptyset$, then $\chi_{\emptyset, \alpha} = \alpha$, which is
not in general a multiple of $\omega_{\alpha}$. However, in the special case
when $I = \Delta - \{ \alpha \}$, we have
$(X, \Phi_{ I \cup \{ \alpha \} }, X^{\vee}, \Phi^{\vee}_{ I \cup \{ \alpha \} })
= (X,\Phi, X^{\vee}, \Phi^{\vee})$, hence 
$\omega_{I\cup \{ \alpha \}, \alpha} = \omega_{\alpha}$, consequently
$\chi_{I,\alpha}$ is a scalar multiple of
$\omega_{\alpha}$ in this case.

\stm\label{coefficients of chi alpha}
{\bf Coefficients of $\chi_{I,\alpha}$.}
Applying co-roots $\gamma^{\vee}$ for $\gamma \in \Delta -I$
and co-weights $\omega^{\vee}_{\beta}$ for $\beta \in I$ to
$\chi_{I,\alpha} = \alpha - \sum_{\beta \in I}b_{\beta} \beta = 
\sum_{\gamma  \in \Delta -I} c_{\gamma}\omega^{\vee}_{\gamma}$
gives
$\langle \chi_{I,\alpha} , \gamma^{\vee}\rangle = c_{\gamma}$, and
$\langle \chi_{I,\alpha} , \omega^{\vee}_{\beta}\rangle = - b_{\beta}$.

\stm\label{basic positive virtual co-character chi alpha}
{\bf Basic positive virtual co-characters $\chi^{\vee}_{I,\alpha}$.}
For $I\subset \Delta$ and $\alpha \in \Delta -I$, we define
the element $\chi^{\vee}_{I,\alpha} \in \Q\otimes X^{\vee}$
as the basic positive virtual co-character corresponding to $(I,\alpha)$
for the dual root datum $(X^{\vee}, \Phi^{\vee}, X, \Phi)$. The duals of the
various statements about $\chi_{I,\alpha}$ hold for the
element $\chi^{\vee}_{I,\alpha}$. In particular, we can write
$$\chi^{\vee}_{I, \alpha} = \alpha^{\vee} - \sum_{\beta \in I} b^{\vee}_{\beta}\beta^{\vee}
= \sum_{\gamma \in \Delta -I}
c^{\vee}_{\gamma} \omega^{\vee}_{\gamma}$$
where each $b^{\vee}_{\beta} \le 0$ for $\beta \in I$ by the dual of
Lemma \ref{R}. Moreover, by the dual of 
Lemma \ref{chi I alpha as a positive multiple}, we have
$$\chi^{\vee}_{I,\alpha} =  \langle \alpha, \chi^{\vee}_{I,\alpha} \rangle 
\omega_{I\cup \{ \alpha \}, \alpha}
\mbox{ and } \langle \alpha, \chi^{\vee}_{I,\alpha} \rangle > 0.$$

\begin{lemma}\label{boundedness and finiteness}
Let  $A \subset \ov{C} \cap (1/n_0) \Q \otimes X^{\vee}$
where $n_0 >0$ is a positive integer. 
Suppose that there exists 
$\gamma \in \Q\otimes \Delta^{\perp}$, an integer
$c$ and a subset $I\subset \Delta$ such that the following three conditions are 
satisfied by each $\mu \in A$.

(1) $\langle \chi_{I,\alpha}, \mu \rangle \le c$ for all $\alpha \in \Delta - I$, 
where $\chi_{I,\alpha}$ are    
the basic positive virtual character associated to $I$ as in 
Lemma \ref{characters positive combination of simple roots}. 

(2) $\langle \beta, \mu \rangle = 0$ for all $\beta \in I$.

(3) $\langle \theta, \mu \rangle = \langle \theta, \gamma \rangle$
for all $\theta \in (\Delta^{\vee})^{\perp}$.

Then the set $A$ is finite.
\end{lemma}

{\bf Proof.} Recall from 
(\ref{basic positive virtual characters for I}) that any $\alpha \in \Delta - I$ can
be expressed as 
$\alpha = \chi_{I,\alpha} + \sum_{\beta \in I}b_{\beta}\beta$ for some $b_{\beta}\in \Q$.
Hence the assumptions (1) and (2) imply the following: 
(4) $\langle \alpha, \mu \rangle \le c$ for all $\alpha \in \Delta - I$.
As $\mu \in \ov{C}$, we also have the following:  
(5) $\langle \alpha, \mu \rangle \ge 0$. 
From (2), (3), (4) and (5), it follows that $S$ is a bounded subset of
$\Q\otimes X^{\vee}$. Now the result follows because
any bounded subset of $(1/n_0) X^{\vee}$ is necessarily finite.
\hfill$\square$

\section{Reductive groups, parabolics, root data}

Proofs of the basic features of the relationship between
root data and split reductive groups (principally
due to Chevalley and Borel) that we will use
are available e.g. in Milne's book [M], which is the reference that
we found most useful for our purpose.

\stm{\bf A split reductive group and its Weyl group.}
Let $k$ be an arbitrary field. A torus over $k$ will mean a group scheme
$T$ such that there exists a field extension $K/k$ such that the base change
$T_K$ is isomorphic over $K$ to $\GG_{m,K}^n$ for some $n\ge 0$
(where $\GG_{m,K}$ is the multiplicative group $\spec K[t,t^{-1}]$).
A split torus over $k$ is a torus $T$ which is isomorphic over $k$ itself
to $\GG_{m,k}^n$ for some $n\ge 0$. A maximal torus $T\subset G$ in a group scheme
$G$ over $k$ is a closed subgroup scheme $T$ which is a torus, such that
$T$ is not properly contained in another such torus $T'\subset G$.

A reductive group over $k$ is a finite-type connected reduced group scheme $G$ over $k$
such that its unipotent radical $R_u(G)$ (which is the unique
maximal connected reduced
closed normal unipotent subgroup scheme of $G$) is trivial. 
A split reductive group over $k$ is a pair
$(G,T)$ where $G$ is a reductive group scheme over $k$, and
$T\subset G$ is closed subgroup scheme which is a split torus,
such that $T$ is a maximal torus in $G$.

The Weyl group $W(G,T)$ of the pair $(G,T)$ is the quotient 
$W(G,T) = N_G(T)/T$ where $N_G(T)$ is the normalizer
of $T$ in $G$. It can be shown that $N_G(T)$ is a closed reduced subgroup scheme of $G$,
and the quotient group scheme $W(G,T)$ is a finite constant group scheme over $k$,
which we can therefore regard as an abstract group. Note that $W(G,T)$ acts on
$T$ by conjugation, and this action is faithful, that is, the homomorphism
$W(G,T) \to Aut_k(T)$ is injective.

\stm{\bf The root datum for a split reductive group.}
To any split reductive group $(G,T)$ over $k$,
there is functorially associated a root datum
$(X^*(T), \Phi, X_*(T), \Phi^{\vee})$ which we now recall. 
Here $X^*(T)$ is the  group 
of all characters (homomorphisms)
$\chi: T\to \GG_m$, and $X_*(T)$ is the group 
of all homomorphisms (one parameter subgroups) $\sigma: \GG_m \to T$.
Composition gives a map $X^*(T)\times X_*(T) \to Hom(\GG_m,\GG_m): (\chi, \sigma)
\mapsto \chi\circ \sigma$, which when composed with the isomorphism
$Hom(\GG_m,\GG_m)\to \Z : (t\mapsto t^n) \mapsto n$ gives the duality pairing
$X^*(T)\times X_*(T) \to \Z$ which is a part of the desired root datum. 

There exists a finite (possibly empty) subset $\Phi\subset X^*(T)$, with
$0\not\in \Phi$, such that
the adjoint action of $T$ on the Lie algebra $\ggg$ of $G$
induces a decomposition $\ggg  = \ggg_0 \oplus (\oplus_{\alpha \in \Phi} \ggg_{\alpha})$
as $T$-representations, 
such that the action of $T$ on the vector space $\ggg_0$ is trivial, while each
$\ggg_{\alpha}$ is a $1$-dimensional vector space on which $T$ acts by multiplication
by the character $\alpha$. This set $\Phi$ is the set of all roots for the
desired root datum. 

The co-roots $\alpha^{\vee} \in X_*(T)$ have the following description.
For each $\alpha \in \Phi$, let $T_{\alpha} = \ker(\alpha : T \to \GG_m)$,   
and let $G_{\alpha} = C_G(T_{\alpha})$ denote the closed subgroup scheme of
$G$ which is the centralizer of $T_{\alpha}$. Then $G_{\alpha}$
is reduced, and it is in fact
a reductive group with maximal torus $T\subset G_{\alpha}$.
Its Weyl group $W(G_{\alpha}, T) = N_{G_{\alpha}}(T)/T$ is of order $2$.
Let $s_{\alpha} \in W(G_{\alpha}, T)$ be its non-identity element. Via its
action on $T$, it defines an automorphism of order $2$ of the group $X^*(T)$,
which sends $\alpha \mapsto - \alpha$.
Hence there exists a unique element $\alpha^{\vee} \in X_*(T)$ such that 
for any $\chi \in X^*(T)$, 
$$s_{\alpha}(\chi) = \chi - \langle \chi, \alpha^{\vee}\rangle \alpha.$$
The map $\Phi \to X_*(T) : \alpha \mapsto \alpha^{\vee}$ can be shown to
be injective. The set of co-roots $\Phi^{\vee} \subset X_*(T)$ is defined to be its image.
This gives a bijection $\Phi \to \Phi^{\vee} : \alpha \mapsto \alpha^{\vee}$, which
is a part of the desired root datum.

With the above definitions, it is a fundamental fact that 
$(X^*(T), \Phi, X_*(T), \Phi^{\vee})$ is indeed a root datum,
naturally associated to the pair $(G, T)$.

In what follows, $(G,T)$ will denote a split reductive group scheme
over a field $k$.

\stm{\bf Characters and virtual characters.} If $H$ is any algebraic group over $k$,
then $\wh{H}$ will denote the abelian group of all characters ($k$-homomorphisms)
$H \to \GG_m$. The elements of $\Q\otimes \wh{H}$ will be called as the virtual characters
of $H$.

\stm\label{Interpretation of global characters and central co-characters}
{\bf $\Q\otimes \wh{G}|_T= (\Phi^{\vee})^{\perp}$ and
$\Q\otimes X_*(Z_0(G)) = \Phi^{\perp}$.} 
For a split reductive $(G,T)$, 
the restriction map $\wh{G} \to X^*(T) : \chi \mapsto \chi|_T$ is an injective
homomorphism. We will denote its image by $\wh{G}|_T \subset X^*(T)$.
The elements of $\wh{G}|_T$ vanish on all co-roots $\alpha^{\vee}$, in fact
$$\Q\otimes \wh{G}|_T = \cap_{\alpha \in \Phi}
\ker(\alpha^{\vee} : \Q\otimes X^*(T) \to \Q)
= (\Phi^{\vee})^{\perp}.$$
Each $\alpha \in \Phi$ is trivial on the center $Z(G)$ of $G$ (which is a
closed subgroup scheme of $T$, not necessarily reduced), and in fact, 
$Z(G) = \cap_{\alpha\in \Phi} \ker(\alpha : T \to \GG_m)$. Note that the connected
component $Z_0(G)$ of the center of $G$ is is a closed
subgroup scheme (not necessarily reduced) of finite index in $Z(G)$. 
Hence
$$\Q\otimes X_*(Z_0(G)) = \Phi^{\perp}.$$

\stm{\bf Semi-simplicity.} A split reductive group $G$ is semi-simple
if and only if $Z_0(G)$ is trivial, equivalently, if and only if $\wh{G}$ is trivial.
It follows from (\ref{Interpretation of global characters and central co-characters}) 
that $G$ is semisimple if and only if its root datum
is semisimple as defined in (\ref{Semi-simplicity.}) above.

\stm{\bf Bases $\Delta \subset \Phi$ correspond to Borels $B \supset T$.}
To any maximal connected solvable reduced closed subgroup scheme
$B$ of $G$ (a Borel subgroup)
with $T\subset B$, we associate the subset $\Phi^+ \subset \Phi$ consisting of all
$\alpha$ such that $\ggg_{\alpha} \subset \bbb$ where $\bbb \subset \ggg$
is the Lie algebra of $B$. Then there exists a unique base $\Delta$ such that
$\Delta \subset \Phi^+$. This defines a  
bijection between the set of all Borels that contain $T$
the set of all bases. In terms of any base $\Delta$, the
identifications  
(\ref{Interpretation of global characters and central co-characters})
take the form
$$\Q \otimes \wh{G}|_T = (\Delta^{\vee})^{\perp} \mbox{ and }
\Q \otimes X_*(Z_0(G)) = \Delta^{\perp}.$$
Equivalently, we have
$$\Q\Delta = X_*(Z_0(G))^{\perp}  \mbox{ and }
\Q\Delta^{\vee} = (\wh{G}|_T)^{\perp}.$$

\stm{\bf Subsets $I\subset \Delta$ correspond to standard parabolics $P\supset B$.}
To any reduced closed subgroup scheme $P\subset G$ with $B\subset P$ (a standard
parabolic) we associate the subset $I_P \subset \Delta$ which consists of
all $\alpha \in \Delta$ such that $\ggg_{-\alpha} \subset \ppp$ where
$\ppp \subset \ggg$ is the Lie algebra of $P$. 
This defines a  
bijection between the set of all standard parabolics and the power set
of $\Delta$. Note that the parabolic $P = B$ corresponds to the empty set
$I_B = \emptyset$, while the parabolic $P = G$ corresponds to $I_G = \Delta$.
The set $I_P$ is called as the set of inverted simple roots of $P$. 
The inverse of the above bijection will be denoted by
$I\mapsto P_I$. \\
{\bf Note.} Any standard parabolic $P$ is automatically connected, and it is its
own normalizer in $G$.

\stm{\bf Characters and virtual characters for $P$.} Let $P\supset B$ be a standard
parabolic, corresponding to the set of inverted roots $I\subset \Delta$.
The restriction map $\wh{P} \to X^*(T)$, under which a character
$\chi : P \to \GG_m$ maps to $\chi|_T \in X^*(T)$, is an injective homomorphism.
We denote its image by $\wh{P}|_T\subset X^*(T)$, or sometimes simply by
$\wh{P} \subset X^*(T)$. We clearly have $\wh{G}|_T \subset \wh{P}|_T\subset X^*(T)$.
The elements of $\Q \otimes  \wh{P}|_T\subset \Q \otimes  X^*(T)$ will be called
as the virtual characters of $P$. For example, if $\beta \in \Delta - I$,
then the fundamental weight $\omega_{\beta}$ is a virtual characters of $P$.
In terms of the root datum, we have
$$\Q \otimes  \wh{G}|_T = (\Delta^{\vee})^{\perp} \mbox{ and } 
\Q \otimes  \wh{P}|_T = (I^{\vee})^{\perp}$$
and so from the decompositions (\ref{Direct sum decomposition-III.}),
we get the following internal direct sum
decompositions.

\stm\label{Decompositions of virtual characters for a parabolic}
{\bf Decomposition of $\Q\otimes X^*(T)$ corresponding to a parabolic.}
With the above identifications, the dual sets of decompositions in 
(\ref{Direct sum decomposition-III.}) take the following form.
\begin{eqnarray*}
X_*(Z_0(G))^{\perp} & = & \Q I_P \oplus (\oplus_{\beta \in \Delta - I_P}\Q\omega_{\beta}),\\
  \Q\otimes \wh{P}|_T & = &
                          (\oplus_{\beta \in \Delta - I_P}\Q\omega_{\beta})
                          \oplus (\Q\otimes \wh{G}|_T),\\
  \Q\otimes X^*(T) & = & \Q \Delta \oplus (\Q \otimes \wh{G}|_T)\\
& = &                    \Q I_P \oplus (\Q\otimes \wh{P}|_T)\\
& = &                       \Q I_P \oplus
                       (\oplus_{\beta \in \Delta - I_P}\Q\omega_{\beta})
                       \oplus (\Q\otimes \wh{G}|_T).\\
  \mbox{Dually},~~~~~~~~~ & & \\
  (\wh{G}|_T)^{\perp} & = & \Q I^{\vee}_P
                        \oplus (\oplus_{\beta \in \Delta - I_P}\Q\omega^{\vee}_{\beta}),\\  
\Q\otimes  X_*(Z_0(P)) & = & (\oplus_{\beta \in \Delta - I_P}\Q\omega^{\vee}_{\beta})
                             \oplus (\Q\otimes X_*(Z_0(G))),\\
  \Q\otimes X_*(T) & = & \Q \Delta^{\vee} \oplus (\Q\otimes X_*(Z_0(G)))\\
                 &= &        \Q I^{\vee}_P \oplus (\Q\otimes X_*(Z_0(P)))\\
                 & = &    \Q I^{\vee}_P \oplus
                         (\oplus_{\beta \in \Delta - I_P}\Q\omega^{\vee}_{\beta})
                         \oplus (\Q\otimes X_*(Z_0(G))).
\end{eqnarray*}                         

\stm\label{character on parabolic trivial on center of G}
{\bf Virtual characters on $P$ that are trivial on $Z_0(G)$.}
Consequently (as an analog of 
\ref{in roots terms character on parabolic trivial on center of G}),
the vector space of virtual characters on $P$ that are trivial on $Z_0(G)$
is given by 
$$
\oplus_{\beta \in \Delta - I_P}\Q\omega_{\beta} =
\Q\Delta \cap (\Q\otimes \wh{P}|_T) = 
X_*(Z_0(G))^{\perp}  
\cap (\Q\otimes \wh{P}|_T).$$
Dually, the vector space  of virtual co-characters of $Z_0(P)$ which are
annihilated by $\wh{G}$ is given by
$$\oplus_{\beta \in \Delta - I_P}\Q\omega^{\vee}_{\beta} =
\Q \Delta^{\vee} \cap (\Q\otimes X_*(Z_0(P))) = 
(\wh{G}|_T)^{\perp} 
\cap (\Q\otimes X_*(Z_0(P))).$$

\stm\label{definition of strongly dominant character}
A {\bf dominant character} for $P$ (which we take in the stronger
sense of Definition 1.14 in [Ram]) is a homomorphism
$\chi : P \to \GG_m$ such that the following two conditions are satisfied:\\
(i) For any simple root $\alpha \in \Delta$ we have
$\langle \chi|_T, \alpha^{\vee}\rangle \ge 0$, where
$\alpha^{\vee} \in X_*(T)$ is the $1$-parameter subgroup of $T$ that is the
co-root associated with $\alpha$, and \\
(ii) The character $\chi$ is trivial on the connected component $Z_0(G)$
of the center $Z(G)$.\\
Equivalently, a character
$\chi : P \to \GG_m$ is dominant in this sense if and only if
the element $\chi|_T \in \Q \otimes X^*(T)$
is a linear combination $\sum_{\alpha \in \Delta - I_P}a_{\alpha}\omega_{\alpha}$
where $a_{\alpha} \in \Q_{\ge 0}$, and $\omega_{\alpha}$ is the fundamental
dominant weight corresponding to any non-inverted root $\alpha$ for $P$.

{\bf Caution about terminology.} The above extra
requirement that a dominant character
$\chi : P\to \GG_m$ should be trivial on $Z_0(G)\subset P$ is special to the
subject of principal bundles. If $G$ is semisimple, this requirement is
automatically satisfied as then $Z_0(G)$ is trivial. However, this  
requirement is not generally a part of the definition of a dominant character
in texts on linear algebraic groups.
Also note that if $\chi : P \to \GG_m$ is a homomorphism, then
the line bundle $L_{-\chi}$ on $G/P$, associated to the
principal $P$-bundle $G\to G/P$ by $-\chi$, is an ample line bundle
if and only if for each simple root $\alpha \in \Delta$ we have
$\langle \chi|_T, \alpha^{\vee}\rangle >0$. Thus, a character $\chi$ on $P$
is dominant in the more customary sense if and only if the line bundle $L_{-\chi}$
lies in the closure of the rationalized ample cone of $G/P$, regardless of whether
$\chi$ is trivial on $Z_0(G)$.

\stm\label{definition of V_P}
{\bf The representation $V_P$.} Let $\omega_{G/P}^{-1} = \det(T_{G/P})$
denote the anti-canonical bundle on
the projective variety $G/P$. This line bundle is very ample, and
gives a closed embedding $G/P \hra {\bf P}(V_P)$ of $G/P$ into the projective space
of lines in $V_P = H^0(G/P ,\omega_{G/P}^{-1})^{\vee}$,
with $\OO_{{\bf P}(V_P)}(-1)|_{G/P} =  \omega_{G/P}$.  
The $G$-action on $G/P$ by left translation has a natural lift to 
$\omega_{G/P}^{-1}$, which makes $V_P$ a finite dimensional
$G$-representation, and makes the above embedding $G$-equivariant.

{\bf Caution}: The $G$-representation $V_P$ is not necessarily irreducible.
For example, let $G = SL(W)$ where $W$ is a $2$-dimensional vector
space over $k$ and $P$ the Borel subgroup that preserves a given
$1$-dimensional subspace $J\subset W$.
Then we have $G/P = {\bf P}(W)$ (the $1$-dimensional projective space of lines
in $W$), $\omega_{G/P}^{-1} = \OO_{{\bf P}(W)}(2)$,
and so $H^0(G/P ,\omega_{G/P}^{-1}) = \sym^2(W^*)$. This (and therefore its dual $V_P$)
is an irreducible representation of $SL(W)$ if and only if $char(k) \ne 2$. 

In any case, $V_P$ has 
a $P$-invariant $1$-dimensional subspace  $J\subset V_P$ 
which corresponds to the $P$-invariant point $eP \in G/P \subset {\bf P}(V_P)$,
where $J$ is just the fiber of $\omega_{G/P}$ over $eP$.
We will denote by $\lambda_P :P \to \GG_m$ the
character by which $P$ acts on $J$.  
As $Z_0(G)\subset P$ the action of $Z_0(G)$ on $V_P$ is trivial,
therefore $\lambda_P$ is dominant in our stronger sense. 
It defines an irreducible sub-$G$-representation $V_{\lambda_P}$ of $V_P$,
and so we get a composite embedding
$$G/P \hra {\bf P}(V_{\lambda_P}) \hra {\bf P}(V_P)$$
under which
$$\OO_{{\bf P}(V_P)}(-1)|_{{\bf P}(V_{\lambda_P})} =
\OO_{{\bf P}(V_{\lambda_P})}(-1) \mbox{ and }
\OO_{{\bf P}(V_{\lambda_P})}(-1)|_{G/P} = \omega_{G/P}.$$

\stm\label{Levi quotient.}
{\bf Levi quotients.} If $B\subset P_I \subset G$ is a standard parabolic,
we will denote by $\phi_I : P_I \to H_I = P_I/R_u(P_I)$ the quotient map, where
$R_u(P_I)$ is the unipotent radical of $P_I$. Note that
$R_u(P) \subset B$. The homomorphism $\phi_I$
(or just the group $H_I$) is known as the Levi
quotient associated to $P$. The morphism $\phi_I$ gives
a closed embedding of $T$ into the quotient 
group $H_I$ (which is reductive), which makes $T$ a maximal torus in $H_I$
and $B_I = B/R_u(P)$ a Borel in $H_I$ containing $T$.
Thus, $(H_I, T)$ is split reductive. Its root datum 
$(X^*(T), \Phi_I, X_*(T), \Phi^{\vee}_I)$ is the root datum associated to
the subset $I \subset \Delta$, defined in \ref{root datum associated to a subset}.
The base for roots of $H_I$ defined by $T\subset B_I \subset H_I$ is the
set $I$. We will denote the corresponding fundamental weights for $H_I$ by
$\omega_{I,\alpha} \in \Q\otimes X^*(T)$ where $\alpha \in I$. {\bf Caution:} In general,
the elements $\omega_{I,\alpha}$ are quite distinct from the fundamental weights
$\omega_{\alpha}$, and these are not even multiples of each other in general.

\stm\label{characters of parabolic and Levi.}
{\bf Identification $\wh{P} = \wh{H}$.} As the kernel $R_u(P)$ of a Levi quotient 
$\phi : P \to H$ is unipotent, any homomorphism $R_u(P) \to \GG_m$ is trivial.
This gives an identifications $\wh{P} = \wh{H}$, and $\Q\otimes \wh{P}|_T = 
\Q\otimes \wh{H}|_T$.

\begin{lemma}\label{the basic positive character generates}
Let $\alpha \in \Delta - I$ and let $\phi : P_{I\cup \{ \alpha \}} \to H_{I\cup \{ \alpha \}}$
be the Levi quotient. Let $\phi(P_I) =
Q' \subset H_{I\cup \{ \alpha \}}$ be the maximal parabolic in $H_{I\cup \{ \alpha \}}$
with $\alpha$ as its only non-inverted root.
With notation as in \ref{basic positive virtual characters for I}, let $n$ be
any positive integer such that $n\chi_{I, \alpha}\in \wh{P}$. Then
the homomorphism $n\chi_{I, \alpha} : P_I \to \GG_m$ factors uniquely
through $\phi : P_I \to Q'$ to define a character $\sigma : Q' \to \GG_m$.
The character $\sigma$ is non-trivial, and its restriction to
$\phi(T) = T$ is a non-negative linear combination of simple
roots for $H_{I\cup \{ \alpha \}}$, that is, of elements of $I \cup \{ \alpha \}$.
Given any character $\tau : Q' \to \GG_m$ whose restriction to
$\phi(T)$ is a non-negative linear combination of elements of
$I \cup \{ \alpha \}$, there exists
a unique rational number $a\ge 0$ such that $\tau|_{\phi(T)} = a \chi_{I, \alpha}$.
\end{lemma}

\proof As $\phi : P_I \to Q'$ is a quotient with kernel the unipotent
group $R_u(P_{I\cup \{ \alpha \} })$, any character $P_I \to \GG_m$ factors uniquely
through $\phi : P_I \to Q'$. The second statement is immediate from
the definition of $\chi_{I, \alpha}$ and the Lemma \ref{R}.
For the final statement, let $\omega_{I\cup \{ \alpha \}, \alpha}$ be the fundamental
weight for $H$ corresponding to the simple root $\alpha$ for $H$. 
Then $\tau|_{\phi(T)}$ is a positive multiple of 
$\omega_{I\cup \{ \alpha \}, \alpha}$, hence also that of $\chi_{I, \alpha}$ by 
(\ref{chi I alpha as a positive multiple}).
\hfill $\square$

\section{Principal bundles and parabolic reductions}

Throughout this section, $X$ will be a smooth irreducible
projective variety over a field $k$, together with a chosen very ample
line bundle $\OO_X(1)$. As before, $(G,T)$ will be a split reductive group over
$k$, together with a choice of a Borel $B$ with $T\subset B \subset G$. Principal
$G$-bundles on $k$-schemes will be assumed to be locally trivial in \'etale
topology (but not necessarily in Zariski topology).

\stm 
Let $E$ be a principal $G$-bundle over a big open subscheme
$W \subset X$, where `big' signifies that
$\dim(X - W) \le \dim(X) -2$, where $\dim(\emptyset) =-\infty$.
A {\bf rational reduction}
$(P,\sigma)$ of the structure group $G$ to a standard parabolic 
$P$ is a section $\sigma: U \to E/P$ 
of $E/P \to X$ over a big open subscheme $U\subset W \subset X$.
The section $\sigma$ has a unique maximal extension to a section
$\sigma_{max}: U_{max}\to E/P$, and we say that two rational reductions
$(P_1,\sigma_1)$ and $(P_2, \sigma_2)$ are equal if $P_1 = P_2$ and
the maximal extension of $\sigma_1$ equals that of $\sigma_2$ (equivalently,
the restrictions of $\sigma_1$ and $\sigma_2$ to the intersection of their domains
are identical to each other, or even more simply, their restrictions to the generic
point of $X$ are identical to each other).

\stm\label{definition of semistability} 
Recall that a principal $G$-bundle $E$ defined on a big open subscheme $W\subset X$
is said to be {\bf semistable} w.r.t. the choice
of $\OO_X(1)$ if for any standard
parabolic $B\subset P\subset G$, any rational reduction
$(P, \sigma: U  \to E/P)$ of $E$
to $P$ and any dominant character
$\chi : P \to \GG_{m,K}$ in the sense of
(\ref{definition of strongly dominant character}), we have 
$$\deg(\chi_*\sigma^*E) \le 0$$
where $\sigma^*E$ is the principal $P$-bundle on $U$ 
defined by the reduction $\sigma$,
and $\chi_*\sigma^*E$ is the $\GG_m$-bundle obtained by
extending its structure group via $\chi: P \to \GG_m$,
which is equivalent to a line bundle on $U$ under the action of $\GG_m$ on
a $1$-dimensional vector space by scalar multiplication. This line
bundle extends uniquely (up to a unique isomorphism)
to a line bundle $L_{\sigma, \chi}$
on $X$, and 
$\deg(\chi_*\sigma^*E)$ is by definition the
degree of $L_{\sigma, \chi}$ w.r.t. $\OO_X(1)$.

\stm\label{diagonalizable implies semistable}
{\bf Case where $G$ is diagonalizable.}
If the connected group $G$ is diagonalizable (that is, $G =T$)
then $Z_0(G) = Z(G) = G$, and the set of roots is empty. The 
only parabolic is $P= G$, for which the only dominant character $\chi$ in the
sense of (\ref{definition of strongly dominant character})
is the trivial character. Thus, the line bundle $\chi_*E$
is trivial so has degree $0$, showing that $E$ is automatically semistable.
This example clearly shows the need for asking for triviality of $\chi$
on $Z_0(P)$ in the definition (\ref{definition of strongly dominant character}) 
of a dominant character for the purpose of defining semistability
of principal bundles.

\stm{\bf Globality assumption for $E$.} From now onwards, we will assume that our
principal $G$-bundle $E$ is defined globally on $X$, unless otherwise indicated
(notably, the principal $P$-bundle $\sigma^*E$
for a rational reduction $(P, \sigma)$, or its
associated principal Levi group bundle, will typically
be defined only on the domain of $\sigma$, not on all of $X$).

\stm\label{Alternative description of reduction}
{\bf Alternative description of $L_{\sigma, \chi}$.} Let
$\sigma: U \to E/P$ be a rational reduction and $\chi : P \to \GG_m$
a dominant character on $P$. 
For later use, we give an alternative description of the line bundle
$L_{\sigma, \chi}$ of (\ref{definition of semistability}),
and show that it comes with a natural injective homomorphism
of sheaves of $\OO_X$-modules $i : L_{\sigma, \chi}\to E(V_{\chi})$, where
$V_{\chi}$ is the irreducible linear representation of $G$ defined
by the dominant weight $\chi$
and $E(V_{\chi})$ is the associated vector bundle to $E$.
A highest weight vector in $V_{\chi}$ (which is unique up to
a non-zero scalar multiple) defines a $1$-dimensional subspace $J_{\chi} \subset
V_{\chi}$. Note that $J_{\chi}$ is invariant under the action of $P\subset G$,
which acts on it by the multiplicative character $\chi$.
The sub-$P$-representation $J_{\chi} \subset V_{\chi}$ defines a line subbundle
$(\sigma^*E)(J_{\chi})$ of the vector bundle
$(\sigma^*E)(V_{\chi}) = E(V_{\chi})|_U$ over $U$. This line bundle
is naturally isomorphic with the line bundle $\chi_*\sigma^*E$ which has
the (uniquely unique) prolongation $L_{\sigma, \chi}$ to $X$, and so 
the inclusion $(\sigma^*E)(J_{\chi})\subset E(V_{\chi})|_U$ uniquely prolongs
to an injective homomorphism of sheaves $i : L_{\sigma, \chi}\to E(V_{\chi})$
which is of rank $1$ at all points of $U$.

\stm\label{recovering the reduction from the homomorphism i}
{\bf Recovering the reduction from the homomorphism
  $i : L_{\sigma,\lambda_P} \to E(V_P)$.}
In what follows, we mainly use a particular choice of $\chi$ for each $P$, namely  
the dominant character $\lambda_P$ on $P$
corresponding to the representation $V_P = H^0(G/P, \omega_{G/P}^{-1})^{\vee}$
of $G$ defined in (\ref{definition of V_P}) in the previous section.
As we have explained, the irreducible representation $V_{\lambda_P}$
with dominant weight $\lambda_P$ is a sub-$G$-representation of $V_P$.
Recall that
for the principal $P$-bundle $G\to G/P$, the associated line bundle $\L$ on
$G/P$ corresponding to $\lambda_P : P \to \GG_m$ has the property that $\L^{-1}$
is very ample, naturally isomorphic as a $G$-homogeneous line bundle to the anti-canonical
bundle $\omega_{G/P}^{-1}$.
To any rational reduction $(P,\sigma)$ of a $G$-bundle
$E$, we have associated in (\ref{Alternative description of reduction})
homomorphism of sheaves 
$L_{\sigma, \lambda_P} \to E(V_{\lambda_P})$ which is fiberwise injective
on a big open $U\subset X$. Composing with the homomorphism
$E(V_{\lambda_P})\hra E(V_P)$ induced by the inclusion $V_{\lambda_P}\hra V_P$
gives a homomorphism 
$$i : L_{\sigma, \lambda_P} \to E(V_P).$$
It is explained in [G-N-2] how to
recover uniquely the reduction $(P, \sigma)$ from the
the line bundle $L_{\sigma, \lambda_P}$ and the homomorphism $i$.
Briefly, the very ample line bundle $\omega_{G/P}^{-1}$ gives an embedding
$G/P \hra {\bf P}(V_P)$ (the projective space of
lines in $V_P$) and hence an embedding $E/P \hra {\bf P}E(V_P)$ over $X$,
and $i$ defines a section of ${\bf P}E(V_P)$ over $U$ which factors via $E/P$,
recovering $\sigma : U \to E/P$. 
In particular, the maximal domain for $\sigma$ is the
open subscheme of $X$ where $i$ has rank $1$ on the fiber at each point.

\stm\label{upper bound on degrees}
{\bf Upper bound on degrees.}
Any vector bundle $\E$ on $(X, \OO_X)$ has a maximal destabilizing
subsheaf $\E_1$ in terms of Gieseker semistability (which is the first
step of its Harder-Narasimhan  filtration, in particular $\E_1 =\E$ if $\E$ is
semistable). 
For any invertible subsheaf $\F\subset \E$, we have $\deg(\F) \le \mu(\E_1)
= \deg(\E_1)/\rank(\E_1) \in \Q$.
Applying this to the vector bundle $\E = E(V_{\chi})$ in the preceding discussion
and its invertible subsheaf $\F = L_{\sigma, \chi}\stackrel{i}{\to} E(V_{\chi}) = \E$,
we get 
$$\deg\chi_*\sigma^*E \le c_{E,P, \chi}$$
where the upper bound
$c_{E,P, \chi} = \mu(\E_1)$ depends only on $(E,P,\chi: P \to \GG_m)$, but is independent
of $\sigma$.

\stm\label{virtual characters of P}
{\bf Dominant virtual characters of $P$.} Recall from
(\ref{Decompositions of virtual characters for a parabolic}) that
the vector space
of virtual characters of a standard parabolic $P$ is
$\Q \otimes \wh{P}|_T = 
(\oplus_{\beta \in \Delta -I_P}\omega_{\beta})
\oplus (\Q\otimes \wh{G}|_T)$, where $I = I_P$. 
We will say that $\chi \in \Q \otimes \wh{P}|_T $
is a {\it dominant virtual character of $P$}, if there exists an integer $a>0$
such that $a\chi$ is a dominant character on $P$ in the sense of
definition (\ref{definition of strongly dominant character}).
This shows that $\chi$ is a dominant virtual character of $P$ if and
only if
$$\chi = \sum_{\beta \in \Delta - I} b_{\beta}\omega_{\beta}
\mbox{ where each } b_{\beta} \in \Q_{\ge 0}.$$

\stm{\bf Degree corresponding to a virtual character of $P$.}
Let $H$ be any algebraic group over $k$, and let $F$ a principal 
$H$-bundle on a large open subscheme $W$ of $X$. Let $\chi \in \Q\otimes \wh{H}$.
There exists an actual character $\rho : H \to \GG_m$ and a rational number
$a \ne 0$ such that $\chi = a \otimes \rho$. Then it can be seen that
the integer $a \deg(\rho_*F)$ is independent of the choice of $(a, \rho)$,
which allows us to define
$\deg(\chi_*F)$ to be $a \deg(\rho_*F)$. In particular, if 
$E$ be a principal $G$-bundle on $X$ and  
$(P,\sigma)$ a rational reduction of $E$, then
$\deg(\chi_*\sigma^*E)$ is defined for any virtual character
$\chi \in \Q \otimes \wh{P}$.
Therefore,
in the definition of semistability of $E$, we can replace the expression `dominant
character for $P$' by the expression `dominant virtual character 
for $P$' to get an equivalent definition.

\stm\label{role of maximal parabolics}
{\bf The role of maximal parabolics.}
For defining semistability, it is enough to consider only the rational
reductions
to the maximal standard parabolics. Such parabolics correspond
to simple roots $\alpha$. The parabolic $P_{\alpha}$
that corresponds to $\alpha$ has $\alpha$ as its unique
non-inverted simple root, 
that is, $I_{P_{\alpha}} = \Delta - \{\alpha \}$.
The corresponding dominant weight
$\omega_{\alpha} \in \Q\otimes X^*(T)$ is the unique non-trivial dominant
virtual character on $P_{\alpha}$ up to positive rational multiples. 
The condition defining semistability of $E$ is equivalent to
the requirement that for all simple roots $\alpha$ and rational reductions
$(P_{\alpha}, \sigma)$ of $E$, we have
$\deg (\omega_{\alpha})_* \sigma^* E \le 0$.
In view of the above discussion, this is in turn equivalent to the
requirement that for all $\alpha$ and all rational
reductions $(P_{\alpha}, \sigma)$, we must have
$\deg(\sigma^*E)(J_{N\omega_{\alpha}}) \le 0$, where $N$ is any positive
integer such that $N\omega_{\alpha} \in X_*(T)$.

\stm{\bf The type of a rational reduction} $(P,\sigma)$ of $E$
is the virtual $1$-parameter subgroup (virtual co-character) 
$\mu_{(P,\sigma)}\in \Q\otimes X_*(T)$ of $T$ 
uniquely defined in terms of the decomposition $\Q\otimes X^*(T)=
(\oplus_{\alpha \in I_P}\Q\otimes I_P) \oplus (\Q\otimes \wh{P}|_T)$
and the natural pairing
$\langle ~,~\rangle : X^*(T) \times X_*(T) \to \Z$ by  
\begin{eqnarray*}
\langle \chi, \mu_{(P,\sigma)}\rangle & = & \left\{
\begin{array}{ll}
\deg(\chi_*\sigma^*E) & \mbox{if }  \chi \in \wh{P},  \\
0                     & \mbox{if } \chi \in I_P.
\end{array}
\right.
\end{eqnarray*}

\stm\label{identity reduction and its type}
{\bf The identity reduction $(G,\id)$ of $E$ and its type $\mu_{(G,\id)}(E)$.}
The identity reduction of $E$ corresponds to the identity section
$\id: X \to E/G = X$, and gives $\id^*E = E$. As $I_G = \Delta$,
the type $\mu_{(G,\id)}(E)$ is contained in the
linear subspace $\Delta^{\perp} = \Q\otimes Z_0(G) \subset \ov{C}$, where
$\ov{C}\subset \Q \otimes X_*(T)$
is the closed positive Weyl chamber defined in (\ref{Weyl chamber}).
Note that $\mu_{(G,\id)}(E)$ depends on $E$, but we will just write it as
$\mu_{(G,\id)}$ for brevity when $E$ is given. 

\stm{\bf Incomparibility of types.}
If $(P,\sigma)$ is any rational reduction of $E$ to a standard parabolic $P$,
then note that $\mu_{(P,\sigma)}(E)  \in \Q\Delta^{\vee} + \mu_{(G,\id)}(E)$. 
Hence if $E$ and $F$ are two principal $G$-bundles and $(P,\sigma)$ and $(Q,\tau)$
are respectively their reductions to standard parabolics, then the types $\mu_{(P,\sigma)}(E)$
and $\mu_{(Q,\tau)}(F)$ are not comparable 
under the partial order $\le$ on $\Q \otimes X^*$ unless 
$\mu_{(G,\id)}(E) = \mu_{(G,\id)}(F)$. When $\mu_{(G,\id)}(E) = \mu_{(G,\id)}(F)$,
we have $\mu_{(P,\sigma)}(E) \le \mu_{(Q,\tau)}(F)$ if and only if
$\deg\chi_*\sigma^*(E) \le \deg\chi_*\tau^*(F)$ for all dominant virtual characters $\chi$ of $P$.

\stm\label{the type of a reduction to a maximal parabolic}
{\bf The type of a reduction to a maximal parabolic $P_{\beta}$.}
Let $\beta \in \Delta$, and let $P_{\beta}$ be the maximal parabolic whose
only non-inverted simple root is $\beta$.
Applying to both sides the linear functionals given by all
$\alpha \in \Delta$ and all $\theta \in \wh{G}$ shows  
that the
type $\mu_{(P_{\beta},\sigma)}$ of any rational reduction $(P_{\beta},\sigma)$ of $E$ 
satisfies the equation 
$$\mu_{(P_{\beta},\sigma)} = \mu_{(G,\id)} +
{\deg ({\omega_{\beta}}_*\sigma^*E) \over
\langle \omega_{\beta}, \omega_{\beta}^{\vee}\rangle} \, 
\omega_{\beta}^{\vee}.$$
As $\langle \omega_{\beta}, \omega_{\beta}^{\vee}\rangle > 0$
(see \ref{Comparison of numbers}), it follows that 
$\mu_{(P_{\beta},\sigma)} \in \ov{C}$ if and only if 
$\deg{\omega_{\beta}}_*\sigma^*E \ge 0$.

\stm\label{type of reduction compared with Levi}
{\bf The type the reduction $(H,\id)$ of the Levi bundle.}
Let $(P,\sigma)$ be a reduction of the $G$-bundle $E$ to a standard 
parabolic $P$ corresponding to $I\subset \Delta$.
Let us denote the type of the reduction more elaborately as
$\mu^G_{(P,\sigma)}(E) \in \Q\otimes X_*(T)$.
Let $\phi : P\to H = P/R_u(P)$ be the Levi quotient, and 
let $E_H = \phi_*\sigma^*E$ be the corresponding $H$-bundle on the large
open subset $U$ of $X$ which is the domain of $\sigma$. 
With respect to the root datum $(X^*(T), \Phi_I, X_*(T), \Phi^{\vee}_I)$
with base $I$ for $H$, and characters for $H$, 
described in (\ref{root datum associated to a subset}), 
(\ref{Levi quotient.}) and (\ref{characters of parabolic and Levi.}), let the type
of $E_H$ be denoted by $\mu^H_{(H,\id)}(E_H)\in \Q\otimes X_*(T)$. 
Then as an immediate consequence of the definitions, we have
$$\mu^G_{(P,\sigma)}(E) = \mu^H_{(H,\id)}(E_H)\in \Q\otimes X_*(T).$$

\stm\label{pull back of a reduction}
{\bf Pull back of a reduction.}
The following observation is valid in the algebraic as well as
topological categories, when the terms `group', `subgroup', `space',
`bundle' etc. are appropriately interpreted.
If $F$ is a principal $P$-bundle on a space $X$ where
$P$ is a group, and if $\phi : P \to H$ is a quotient, then the associated bundle
$\phi_*F$ can be identified with the quotient $F/R$ where $R = \ker(\phi)$
(we assume $R$ is closed).
Let $K\subset H$ be a subgroup, let $Q = \phi^{-1}(K) \subset P$,
and let $(K,s)$, where
$s : X \to (\phi_*F)/K$ is a section, be a reduction of the structure group of 
$\phi_*F$ from $H$ to $K$. Then as
$(\phi_*F)/K = (F/R)/(Q/R) = F/Q$, we can regard $s$ as a section $s : X \to F/Q$
which defines a reduction $(Q,s)$ 
of the structure group of $F$ from $P$ to $Q$. We will call $(Q,s)$ as
the pullback of the reduction $(K, s)$.
If $\phi' : Q\to K$
denotes the restriction of $\phi$, then we have a natural identification
$\phi'_* (s^*F) = s^*\phi_*F$.

\stm\label{Composite reduction}
{\bf Composite reduction.} Let $G$ be a group, $Q\subset P \subset G$
subgroups, $E$ a principal $G$-bundle on a space $X$, $\sigma : X \to
E/P$ a reduction of the structure group of $E$ to $P$ and $\tau : X \to (\sigma^*E)/Q$
a reduction of the structure group of $\sigma^*E$ to $Q$. The reduction $\sigma$
gives rise to a $P$-equivariant embedding $i_{\sigma}: \sigma^*E \hra E$.
Hence $\tau$ gives a section $X \to E/Q$ which is the composite
$$X \stackrel{\tau}{\to}  (\sigma^*E)/Q \stackrel{\ov{i_{\sigma}}}{\to} E/Q$$
and so we get a reduction  of the structure group of $E$ to $Q$.
We call this as the composite of the reductions $\sigma$ and $\tau$, and
denote it by $\sigma\circ \tau$ for simplicity of notation.
Note that $(\sigma \circ \tau)^*E = \tau^*\sigma^*E$ as a $Q$-bundle.

\begin{lemma}\label{type of the pullback of reduction to a maximal parabolic of Levi}
Let $(P,\sigma : U\to E/P)$ be a rational reduction of a $G$-bundle $E$
to a standard parabolic
$P$ that corresponds to $I \subset \Delta$. Let $E_H$ be the corresponding
$H$-bundle $\phi_*\sigma^*E$ where $\phi : P \to H$ is the Levi quotient.
Let $\alpha \in I$ and let $Q_{\alpha} \subset H$ be the maximal
parabolic in $H$ whose only non-inverted simple root is $\alpha \in I$ w.r.t
the base $I$ for $H$. Let $(Q_{\alpha},s : W \to E_H/Q_{\alpha})$ be a rational
reduction of $E_H$ to $Q_{\alpha}$ defined over a large open subscheme $W\subset U$.
Let $P' = P_{I - \{ \alpha \}}\subset P$, which is the inverse image of $Q_{\alpha}$
under $\phi : P\to H$. Let $(P',s)$ be the rational reduction of $\sigma^*E$ that is
the
pull back of the rational reduction $(Q_{\alpha}, s)$ of $E_H$
as in (\ref{pull back of a reduction}), and
let $(P', \sigma\circ s)$ be the resulting composite
reduction of $E$ to $P'$ as in (\ref{Composite reduction}).
Then we have
$$\mu^G_{(P', \sigma\circ s) }(E) = \mu^G_{(P,\sigma)}(E) + 
{\deg ((\omega_{I, \alpha})_*s^*E_H) \over
\langle \omega_{I,\alpha}, \omega_{I, \alpha}^{\vee}\rangle}
\omega^{\vee}_{I, \alpha}$$
where $\omega^{\vee}_{I, \alpha}$ denotes the fundamental co-weight for
the simple root $\alpha$ for $H$ as given in
(\ref{root datum associated to a subset}).
\end{lemma}

\proof By (\ref{type of reduction compared with Levi}), we have
$$\mu^H_{(H,\id)}(E_H) = \mu^G_{(P,\sigma)}(E).$$
By (\ref{the type of a reduction to a maximal parabolic}) applied to
the reductive group $H$, its maximal parabolic $Q_{\alpha} \subset H$,
the fundamental weight $\omega_{I, \alpha}$ for $H$ for the simple root
$\alpha\in I$, the $H$-bundle $E_H$
and its rational reduction $(Q_{\alpha}, s : W \to E_H/Q_{\alpha})$, we get
$$\mu^H_{(Q_{\alpha},s)}(E_H) = \mu^H_{(H,\id)}(E_H) +
{\deg ((\omega_{I, \alpha})_*s^*E_H) \over
\langle \omega_{I,\alpha}, \omega_{I, \alpha}^{\vee}\rangle} \, 
\omega_{I, \alpha}^{\vee}.$$
From the above two equations, we get
$$~~~~~~~~~
\mu^H_{(Q_{\alpha},s)}(E_H)  = \mu^G_{(P,\sigma)}(E)+
{\deg ((\omega_{I, \alpha})_*s^*E_H) \over
\langle \omega_{I,\alpha}, \omega_{I, \alpha}^{\vee}\rangle} \, 
\omega_{I, \alpha}^{\vee}
~~~~~~~~~~~\ldots (*)$$
We now apply (\ref{pull back of a reduction}) taking $F = \sigma^*E$,
which is a principal $P$-bundle for $P = P_I$. We take $\phi : P \to H$
to be the Levi quotient and $K = Q_{\alpha}\subset H$, so that
$\phi^{-1}(K) = P_{I -\{ \alpha \}}$. 
We are given a rational reduction $(Q_{\alpha}, s)$ of $E_H =
\phi_* \sigma^*E = \phi_*F$. Hence by (\ref{pull back of a reduction}),
we get a reduction $(P_{I - \{ \alpha \}}, s)$ of $\sigma^*E$, which
is the pullback of the reduction $(Q_{\alpha}, s)$ of $E_H$.
Let $\phi' : P_{I - \{ \alpha \}}\to Q_{\alpha}$ denote the restriction of $\phi$
to $P_{I - \{ \alpha \}} \subset P$. Then by (\ref{pull back of a reduction}),
we have $\phi'_*s^*\sigma^*E = s^*\phi_*\sigma^*E$.

Let $i: P_{I - \{ \alpha \}}\hra P$ be the inclusion. Then $i_*s^* \sigma^*E = \sigma^*E$.
Consider the composite reduction $r = \sigma \circ s$ of $E$ to $P_{I - \{ \alpha \}}$
as defined in (\ref{Composite reduction}), with $r^*E = s^*\sigma^*E$.
Hence $i_*r^*E = i_*s^*\sigma^*E = \sigma^*E$, and
$\phi'_* r^* E = \phi'_*s^*\sigma^*E = s^*\phi_*\sigma^*E$.

Let $\psi: Q_{\alpha} \to H_{\alpha}$ be the Levi quotient of $Q_{\alpha}$.
Let $\psi_*s^*E_H$ be the associated principal $H_{\alpha}$-bundle.
By (\ref{pull back of a reduction}),
$\psi_*s^*E_H = \phi'_*s^*\sigma^*E = \psi'_*r^*E$, which
is the extension of
the principal $P_{I - \{ \alpha \}}$-bundle $r^*E$ under the Levi 
quotient $\psi' = \psi \circ \phi'  : P' \to Q_{\alpha} \to H_{\alpha}$.
By  (\ref{type of reduction compared with Levi}) applied
to the both the Levi quotients $\psi: Q_{\alpha} \to H_{\alpha}$ and
$\psi': P_{I - \{ \alpha \} } \to H_{\alpha}$, we have
$$\mu^H_{(Q_{\alpha},s)}(E_H) = \mu^{H_{\alpha}}_{(H_{\alpha},\id)}(\psi_*s^* E_H)
= \mu^{H_{\alpha}}_{(H_{\alpha},\id)}(\psi'_*r^*E) 
= \mu^G_{(P_{I - \{ \alpha \} },r)}(E)$$
and hence 
$$~~~~~~~~~~~~ \mu^H_{(Q_{\alpha},s)}(E_H) =  \mu^G_{(P_{I - \{ \alpha \} },r)}(E)
~~~~~~~~~~~~~~~\ldots (**)$$
From the above two displayed equations $(*)$ and $(**)$, we have 
$$\mu^G_{(P_{I - \{ \alpha \} },r)}(E) = \mu^G_{(P,\sigma)}(E)+
{\deg ((\omega_{I, \alpha})_*s^*E_H) \over
\langle \omega_{I,\alpha}, \omega_{I, \alpha}^{\vee}\rangle} \, 
\omega_{I, \alpha}^{\vee}$$
which completes the proof.
\hfill$\square$

\stm\label{Maximality of reduction implies semistability of Levi bundle}
{\bf Maximality of reduction type implies semistability of Levi bundle.}
Let $(P,\sigma)$ be a rational reduction of a principal $G$-bundle $E$
to a standard parabolic $P$, such that its type $\mu_{(P,\sigma)}$ is maximal
among types of rational reductions, w.r.t. the partial order
(\ref{definition of partial order}) on $\Q\otimes X_*(T)$ that corresponds to
the choice of $T\subset B\subset G$.
Then the $H$-bundle $E_H = \phi_*\sigma^*E$ is semistable, where $\phi: P\to H$ is
the Levi quotient of $P$. For, if $E_H$ is not semistable, then
by (\ref{role of maximal parabolics}) 
there exists a maximal parabolic $Q_{\alpha} \subset H$ and a
rational reduction $s : W \to E_H/Q_{\alpha}$ which contradicts
semistability of $E_H$, that is,
$\deg (\omega_{I, \alpha})_*s^*E_H >0$. 
By (\ref{positivity of fundamental wts for I w.r.t. Delta}),
we have $\omega_{I, \alpha}^{\vee} >0$ in terms of the partial order on $\Q\otimes X_*(T)$
that corresponds to $T\subset B \subset G$.
Hence by Lemma \ref{type of the pullback of reduction to a maximal parabolic of Levi}
we have a rational reduction $(P_{I - \{ \alpha \} },\sigma\circ s)$ of $E$ 
such that
$$\mu^G_{(P_{I - \{ \alpha \} },\sigma\circ s)}(E) > \mu^G_{(P,\sigma)}(E).$$
This contradicts the assumption that $\mu^G_{(P,\sigma)}(E)$ is
maximal.
\hfill $\square$

\begin{lemma}\label{types lattice}
There exists an integer $n_0 >0$, which is 
independent of the projective variety $X$, the line bundle $\OO_X(1)$
and the principal $G$-bundle $E$ on it, 
such that
for all rational reductions $\sigma$ of $E$ 
to standard parabolics $P$, we have  
$$\mu_{(P,\sigma)} \in {1\over n_0} X_*(T)\subset \Q \otimes X_*(T)$$  
\end{lemma}

\proof Let $\{ \theta_i \,|\, i = 1,\ldots, n\}$ be a $\Z$-basis for
the free abelian group $\wh{G}$. 
Let $N>0$ be an integer such that for any standard parabolic $P$,
we have $N \omega_{\alpha} \in \wh{P}|_T \subset X^*(T)$ for all
$\alpha \in \Delta - I_P$. For any $I\subset \Delta$, the  
vector space $\Q\otimes X^*(T)$ has a $\Q$-linear basis
${\cal B}_P  = I \cup \{ N\omega_{\alpha} \,|\, \alpha \in \Delta - I\}
\cup \{ \theta_i \,|\, i = 1,\ldots, n\}$,
and this basis ${\cal B}_P$ lies in $X^*(T) \subset \Q\otimes X^*(T)$.
Hence by the structure theorem for finitely generated abelian groups,
there exists an integer $n_P$ such that
$$n_P X^*(T) \subset \Z {\cal B}_P \subset X^*(T).$$
Now, if $v \in {\cal B}_P$, then by definition of $\mu_{(P,\sigma)}$
we have $\langle v, \mu_{(P,\sigma)} \rangle \in \Z$.
Hence if $x \in X^*(T)$ then as $n_Px \in \Z{\cal B}_P$ we have 
$$\langle x , n_P \mu_{(P,\sigma)} \rangle =
\langle n_P x, \mu_{(P,\sigma)} \rangle\in \Z.$$
This shows that $n_P\mu_{(P,\sigma)} \in X_*(T)$, which is the dual of $X^*(T)$.
Taking $n_0>0$ to be a common multiple of the $n_P$ as $P$ varies over the
standard parabolics of $G$, the lemma follows.
\hfill$\square$

\stm {\bf Canonical reductions.}
Recall that a {\it canonical reduction} of $E$ is  
a rational reduction $\sigma : U \to E/P$ 
of the structure group of $E$ to a standard parabolic $P\subset G$
for which the following two conditions (C-1) and (C-2) hold:

{\bf (C-1)} If $\phi : P \to H = P/R_u(P)$ is the Levi quotient of $P$ 
then the principal $H$-bundle 
$\phi_*\sigma^*E$ defined on $U$ is semistable. 

{\bf (C-2)} For any non-trivial character $\chi: P \to \GG_m$ 
whose restriction to the chosen maximal torus $T\subset B \subset P$ 
is a linear combination $\sum_{\alpha \in \Delta} n_{\alpha}\alpha$ 
of simple roots where each $n_{\alpha} \ge 0$ and at least
one $n_{\alpha} \ne 0$, we have $\deg(\chi_*\sigma^*E) > 0$.

\stm{\bf The reduction $(G,\id)$ is canonical if and only if $E$ is semistable},
as follows immediately from the definitions and the observation that
the only character $\chi: G \to \GG_m$ 
whose restriction to $T$ 
is a linear combination $\sum n_i\alpha_i$ 
of simple roots $\alpha_i \in \Delta$ is the trivial character $1$
(means $\chi|_T = 0$ in additive notation).

\stm\label{C 2 implies type is in Weyl chamber}
{\bf If (C-2) is satisfied by $(P,\sigma)$ then $\mu_{(P,\sigma)} \in \ov{C}$.}
If $P = G$, then $(P,\sigma) = (G,\id)$,
and we have seen in statement \ref{identity reduction and its type}
that $\mu_{(G,id)} \in \Q\otimes Z_0(G) \subset \ov{C}$. So let $P\ne G$,
and hence the set $I\subset \Delta$
of inverted roots for $P$ is a proper subset of $\Delta$.
For $\alpha \in \Delta -I$, consider the positive virtual character 
$\chi_{\alpha} \in \Q \otimes \wh{P}|_T$ of 
(\ref{basic positive virtual characters for I}), 
which equals $\alpha - \sum_{\beta \in I}b_{\beta}\beta$ where each $b_{\beta} <0$.
As $\alpha = \chi_{\alpha} + \sum_{\beta \in I}b_{\beta}\beta$
and as $\langle \beta, \mu_{(P,\sigma)} \rangle = 0$ for all
$\beta \in I$, we get 
$$\langle \alpha, \mu_{(P,\sigma)} \rangle =
\langle \chi_{\alpha}, \mu_{(P,\sigma)}\rangle > 0$$
where the last inequality follows from (C-2) as
$\chi_{\alpha}>0$ in terms of the partial order
(\ref{definition of partial order}) on $\Q \otimes X^*(T)$. 
This shows that $\mu_{(P,\sigma)} \in \ov{C}$.

\stm {\bf Comparison with vector bundles.} Let $E$ be a principal $GL_n$-bundle
where $n\ge 1$, 
and let $\E$ be the corresponding vector bundle on $X$. Then semistability
of the principal bundle $E$ is equivalent to the concept of $\mu$-semistability
(slope semistability) for $\E$. Let $1\le r_1 < r_2 < \ldots < r_{\ell} = n$,
and let $P\subset GL_n$ be the parabolic that preserves the flag
$0 \subset k^{r_1} \subset \ldots \subset k^{r_{\ell}} = k^n$.
Then any rational $P$-reduction of $E$
corresponds to a filtration of $\E$ by coherent subsheaves
$0 \subset \E_1 \subset \ldots \subset \E_{\ell} = \E$
such that the quotient $\E_{i+1}/\E_i$ is
torsion-free of rank $r_{i+1} - r_i$
for $i = 0,\ldots, \ell -1$. The condition
(C-1) says that each $\E_{i+1}/\E_i$ is $\mu$-semistable. The condition (C-2) says
that we should have strict inequalities
$\deg(\E_i)/r_i > \deg(\E_{i+1})/r_{i+1}$ for $i = 1, \ldots, \ell -1$.
This means a rational $P$-reduction of $E$ is canonical if and only if
the corresponding filtration 
$0 \subset \E_1 \subset \ldots \subset \E_{\ell} = \E$ is the 
Harder-Narasimhan filtration of $\E$ w.r.t. $\mu$-semistability.

\stm{\bf Converse of \ref{C 2 implies type is in Weyl chamber} is false.}
In the above, let $\E = \L \oplus \L$ be a vector bundle of rank $2$, where
$\L$ is a line bundle on $X$. This corresponds to a principal $GL_2$
bundle $E$ on $X$, and the filtration $\E_1 = (\L\oplus 0)  \subset \E$ corresponds to
a reduction $(B,\sigma)$
to the Borel $B\subset GL_2$ of upper triangular matrices. 
The type $\mu_{(B,\sigma)}$ lies in $\ov{C}$. Moreover, the reduction satisfies
(C-1). But this is not a canonical reduction, as (C-2) is not satisfied.
(Note that by the Theorem (\ref{existence and uniqueness}) 
of uniqueness of a canonical reduction, 
as $E$ is already $G$-semistable, $(B,\sigma)$ could not have been canonical.)

\begin{proposition}\label{notherian condition on types}
{\bf The poset of reduction types for a given $E$ is noetherian.} 
Let $E$ be a principal $G$-bundle on $X$ and let $P$
be any standard parabolic. Let $S_P\subset \Q\otimes X_*(T)$
be the set of types $\mu_{(P,\sigma)}$
of all possible rational reductions $(P,\sigma)$
of the structure group of $E$ to $P$. 
Then the poset $S_P$ is noetherian w.r.t the partial order on $\Q\otimes X_*(T)$.
Therefore, the finite union $S = \cup S_P$ over all standard parabolics $P$ 
is again noetherian. In particular, as $S$ is nonempty (as $S_G$ is a singleton),
there exists a maximal element in $S$.
\end{proposition}

\proof Let $P$ correspond to $I\subset \Delta$. Consider the decomposition
$$\Q\otimes X^*(T) = 
\Q I \oplus (\oplus_{\beta \in \Delta - I}\Q\omega_{\beta}) \oplus (\Q\otimes \wh{G}|_T).$$
By definition of $\mu_{(P,\sigma)}$, we have
$\langle I, \mu_{(P,\sigma)} \rangle =0$, while $\langle \gamma, \mu_{(P,\sigma)} \rangle =
\langle \gamma, \mu_{(G,\id)}\rangle$ for all $\gamma \in \Q\otimes \wh{G}|_T$. 
Consider the basis $\chi_{I,\alpha}$ for $\oplus_{\beta \in \Delta - I}\Q\omega_{\beta}$
given by (\ref{basic positive virtual characters for I}), where $\alpha$ varies over $\Delta -I$.
Let $x,y \in \Q \otimes X_*(T)$ such that $\langle I,x\rangle =  
\langle I, y \rangle = 0$ and $\langle \gamma,x\rangle =  
\langle \gamma,y \rangle$ for all $\gamma \in \Q\otimes \wh{G}|_T$.
Then $x< y$ if and only if $\langle \chi_{I,\alpha}, x\rangle \le  
\langle\chi_{I,\alpha}, y \rangle$ for all $\alpha \in \Delta - I$ and 
$\langle \chi_{I,\alpha}, x\rangle <  \langle\chi_{I,\alpha}, y \rangle$
for at least one $\alpha \in \Delta - I$.
Hence if $x_1 < x_2 < \ldots $ is a strictly increasing sequence in $S_P$ then there exists
some $\alpha \in \Delta -I$ such that the sequence $\langle\chi_{I,\alpha}, x_n\rangle$ in $\Q$ has
strictly increasing subsequence.
As there exists an integer $r \ge 1$ such that 
$r \chi_{I,\alpha}$ is an actual character on $P$, we have
$\langle r \chi_{I,\alpha}, \mu_{(P,\sigma)}\rangle = \deg(r \chi_{I,\alpha})_*\sigma^*E 
\in \Z$ for any $(P,\sigma)$, and so
$\langle\chi_{I,\alpha}, x_n\rangle$ has a
strictly increasing subsequence in ${1\over r}\Z \subset \Q$.
Hence $\langle\chi_{I,\alpha}, x_n\rangle$ is
not bounded above. This contradicts the statement (\ref{upper bound on degrees}), completing the proof.
\hfill $\square$

{}

We are now ready to prove the following well-known fundamental fact.
Our proof of existence is new in all dimensions, while 
that of uniqueness is new in dimensions $\ge 2$, as far as we know.
The strategy of our proof has a portion similar to the argument 
in [Bi-Ho].

\begin{theorem}\label{existence and uniqueness}
Let $G$ be a split reductive group over a field $k$ of arbitrary characteristic. 
Then any principal $G$-bundle on a smooth irreducible projective variety
$X$ over $k$ admits at least one canonical reduction. If moreover
$k$ is algebraically closed, or more generally, if the hypothesis ${\bf (*)}$
is satisfied and $X$ is geometrically irreducible,
then the canonical reduction is unique.
\end{theorem}

\proof 
{\it Existence of at least one canonical reduction.}
We put a partial order on the set of all rational reductions
$(P,r)$ of $E$ to standard parabolics $P$ by defining
$(P,r) \le (Q,s)$ if $P\subset Q$ and $\mu_{(P,r)} \le \mu_{(Q,s)}$ in $\Q\otimes X_*(T)$. 
As chains of inclusions
$P_1\subset P_2 \subset \ldots$ become stationary, and as we have already
shown in Lemma \ref{notherian condition on types}
that the types $\mu_{(P,r)}$ form a noetherian poset in $\Q\otimes X_*(T)$,
it follows that there exists a maximal $(P,r)$. We will show that any
such $(P,r)$ is a canonical reduction, thereby proving the existence.

We have already proved in
(\ref{Maximality of reduction implies semistability of Levi bundle}) above
that any $(P,r)$ for which $\mu_{(P,r)}$ is maximal
satisfies the condition {\bf (C-1)} in the definition of a canonical reduction. \\
{\bf Note.} This part does not need $P$ to be maximal w.r.t. inclusion among
all those which maximize $\mu_{(P,r)}$.

Next, we prove that any maximal $(P,r)$ satisfies {\bf (C-2)}. 
Let $P$ correspond to $I\subset \Delta$.
To prove {\bf (C-2)}, by (\ref{characters positive combination of simple roots})
we just need to show that for any $\alpha \in \Delta -I$,
the basic positive virtual character $\chi_{I, \alpha}$ defined in
(\ref{basic positive virtual characters for I})
has the property
that $\deg(\chi_{I, \alpha})_*r^*E > 0$. Let $N$ be a positive integer, such
that $N\chi_{I, \alpha}$ is an actual character $P \to \GG_m$.

Let $P\subset Q$ where the parabolic $Q$ corresponds to 
$I_Q = I_P \cup \{ \alpha \}$, and let $(Q,s)$ be the induced reduction 
of $E$, so that $s^*E = i_*r^*E$ where $i: P\hra Q$ is the inclusion homomorphism.
Let $\phi : Q\to H$ be the Levi quotient of $Q$, so that $\phi(P) =
Q_{\alpha} \subset H$
is the maximal parabolic of $H$ corresponding to the only non-inverted root
$\alpha \in I \cup \{ \alpha \}$ (where recall that
$I \cup \{ \alpha \}$ is the set of simple roots for $H$).
Then by Lemma \ref{chi I alpha as a positive multiple}, we have the equality
$\chi_{I, \alpha} = a \omega_{I \cup \{ \alpha \}}$ where
$a = \langle \chi_{I, \alpha} , \alpha^{\vee} \rangle >0$.
By Lemma \ref{type of the pullback of reduction to a maximal parabolic of Levi},
the types $\mu_{(P,r)}$ and $\mu_{(Q,s)}$ are related by
$$\mu^G_{(P, r) }(E) = \mu^G_{(Q,s)}(E) + 
{\deg ((\omega_{I\cup \{ \alpha \}, \alpha})_*s^*E_H) \over
\langle \omega_{I\cup \{ \alpha \},\alpha}, \omega_{I\cup \{ \alpha \}, \alpha}^{\vee}\rangle}
\omega^{\vee}_{I\cup \{ \alpha \}, \alpha}$$
where $\omega^{\vee}_{I\cup \{ \alpha \}, \alpha}$ denotes the fundamental co-weight for
the simple root $\alpha$ for $H$ as given in
(\ref{root datum associated to a subset}).

Note that $\chi_{I, \alpha}$ and $\omega_{I\cup \{ \alpha \} , \alpha}$ are
virtual characters of both $P$ and $H$ under the identification
(\ref{characters of parabolic and Levi.}) induced by $\phi: P\to H$.
Now suppose we have $\deg(\chi_{I, \alpha})_*r^*E \le 0$.
But note that $\deg(\chi_{I, \alpha})_*r^*E$ is a positive multiple of
$\deg ((\omega_{I\cup \{ \alpha \}, \alpha})_*s^*E_H)$. Hence if
$\deg(\chi_{I, \alpha})_*r^*E \le 0$, then we have 
an inequality
$$\mu_{(P,r)} \le \mu_{(Q,s)}.$$
By maximality of $\mu_{(P,r)}$ this means $\mu_{(P,r)} = \mu_{(Q,s)}$.
But then by maximality of $P$ under inclusions, we must have $P=Q$,
which is a contradiction. Hence for every $\alpha \in \Delta -I$ we must have
$\deg(\chi_{I, \alpha})_*r^*E > 0$.
Hence {\bf (C-2)} holds by Lemma \ref{characters positive combination of simple roots}.

This completes the proof of the existence of at least one canonical reduction.

\bigskip

{\it Proof that under the hypothesis ${\bf (*)}$ there is at most one
  canonical reduction.}

Assume first that $k$ is algebraically closed. 
Let $\OO_X(1)$ be the chosen very ample line bundle on $X$.
Let $(\sigma_1: U_1 \to E/P_1)$ and $(\sigma_2: U_2 \to E/P_2)$ be canonical
reductions of $E$, where $P_1$ and $P_2$ are standard parabolics. Let
$U = U_1\cap U_2 \subset X$, which is a large open subscheme in $X$.
Now consider a complete-intersection curve $C\subset X$
obtained by intersecting smooth hypersurfaces in the linear system on $X$ 
corresponding to a 
power $\OO_X(d)$ of $\OO_X(1)$. 
Note that if $d$ is chosen large enough and $C$ is chosen to be sufficiently general,
then $C$ will be smooth over $k$ and will be entirely contained in the
large open subset $U\subset X$.
The restrictions $\sigma_1|_U : C \to (E|_C)/P_1$ and $\sigma_2|_U : C \to (E|_C)/P_2$
will be parabolic reductions of $E|_C$, which by the semistable restriction
theorem (Theorem 12 in [Gu]) will each define a canonical reduction of $E|_C$.
By the uniqueness of canonical reductions over curves proved by Behrend in [Be],
these reductions coincide on $C$, in particular, the standard parabolics $P_1$
and $P_2$ which underlie the two canonical reductions.
We will write $P_1=P_2 = P$.

Recall from the Remark \ref{recovering the reduction from the homomorphism i} that 
a reduction $(P,\sigma)$ can be uniquely recovered from its associated
homomorphism $i : L_{P,\sigma} \to E(V_P)$.
Let the two distinct canonical reductions $(P,\sigma_i)$, $i=1,2$, of $E$ correspond 
to homomorphisms $i_r: L_i \rightarrow E(V_P)$, $r=1,2$, where
$L_i = L_{P,\sigma_i}$ and $V_P = H^0(G/P, \omega^{-1}_{G/P})^{\vee}$.
The homomorphisms $i_1$ and $i_2$ are fiberwise injective on 
the large open set $U\subset X$.
As the reductions are equal on $C$, there exists an $\OO_C$-linear isomorphism
$\psi : L_1|_C \to L_2|_C$ with $i_2|_C \circ \psi = f_1|_C$. 

When $C$ has sufficiently high degree, by the Enriques-Severi-Zariski lemma,
the restriction homomorphism
$H^0(X,L_2\otimes L_1^{-1}) \to H^0(C,L_2\otimes L_1^{-1}|_C)$ is an isomorphism. 
Hence the section $\psi \in H^0(C,L_2\otimes L_1^{-1}|_C)$
can be uniquely lifted to a section $\tilde{\psi}\in H^0(X,L_2\otimes L_1^{-1})$,
that is, to a homomorphism $\tilde{\psi}: L_1 \rightarrow L_2$ on $X$.
As the vanishing divisor of this homomorphism does not intersect $C$, it follows
that $\tilde{\psi}$ will in fact be an isomorphism over $X$. 
Hence we can without loss of generality assume that the two canonical
reductions over $X$ correspond to homomorphisms $i_1: L \to E(V_P)$ and 
$i_2: L \to E(V_P)$ which share a common line bundle $L$ as their domain.
Hence we have two sections $i_r \in H^0(X, E(V_P)\otimes L^{-1})$, $r=1,2$, 
which have the same restriction to $C$.
Again, if $d$ is sufficiently large, the restriction homomorphism
$H^0(X, E(V_P)\otimes L^{-1})\to H^0(C,E(V_P)\otimes L^{-1}|_C)$ is an isomorphism, 
hence $i_1 = f_2$ as desired. Therefore $(P,\sigma_1) = (P,\sigma_2)$
by Remark  \ref{recovering the reduction from the homomorphism i}.
This completes the proof when $k$ is algebraically close.

Finally, we consider the case of a general $k$. We assume that $X$ is geometrically
irreducible over $k$, and moreover the hypothesis ${\bf (*)}$ is satisfied. 
If $E$ were to admit two different canonical reductions on $X$, then their
pull backs to $X_{\ov{k}}$ under base change to an algebraic closure $\ov{k}$
of $k$ would be distinct
reductions of $E_{\ov{k}}$ by faithful flatness of $\ov{k}/k$. 
By hypothesis ${\bf (*)}$, the pulled back reductions are canonical over $\ov{k}$.
But this would contradict the uniqueness proved for $G_{\ov{k}}$-bundles on
the smooth irreducible projective variety $X_{\ov{k}}/\ov{k}$.
Hence uniqueness holds over $k$, which completes the proof of the theorem. 
\hfill $\square$

\stm{\bf The Harder-Narasimhan type of a principal bundle.}
Henceforth, we will assume that the hypothesis ${\bf (*)}$ is satisfied. 
If $(P,\sigma)$ is the canonical reduction of $E$
(which exists and is unique by Theorem \ref{existence and uniqueness}
under the hypothesis ${\bf (*)}$),
then its type $\tau = \mu_{(P,\sigma)}\in \Q\otimes X_*(T)$
is called the {\it Harder-Narasimhan type of $E$}, or just the {\it type} of $E$.

\begin{corollary}\label{maximality of HN type}
  {\bf Maximality of HN type}
  For a given $E$, 
  let the subset $A \subset \Q\otimes X_*(T)$ consist of
  all $\tau$ such that $ \tau = \mu_{(P,\sigma)}(E)$
for some rational reduction $(P,\sigma)$ of $E$ where $P$ is a standard parabolic.
Let $A$ be given the partial order $\le$ induced from $\Q\otimes X_*(T)$. 
Then the Harder-Narasimhan type $\HN(E)$ of $E$ is 
the unique maximum element of the poset $(A, \le)$.
\end{corollary}  

\proof We first show how to possibly 
enlarge the parabolic $P$ to $P\subset P'$, without changing
the type of the reduction.
Given a reduction $(P,\sigma)$ of $E$, let
$J\subset \Delta$ consist of all $\alpha$ such that
$\langle \alpha, \mu_{(P,\sigma)} \rangle =0$. By definition of
$\mu_{(P,\sigma)}$, we have $I_P \subset J$. Let $P'$ be the standard parabolic
with $I_{P'} = J$, so $P\subset P'$. We get a quotient $E/P \to E/P'$, and its
composition with $\sigma: U \to E/P$ gives a section $\sigma' : U \to E/P'$.
By its construction, the new reduction $(P',\sigma')$ has the same type
$\mu_{(P',\sigma')} =  \mu_{(P,\sigma)}$. 
If $(P,\sigma)$ is a canonical
reduction of $E$, then $(P',\sigma') = (P,\sigma)$ by the maximality
of a canonical reduction under the partial order on the set of
all rational reductions defined in the proof
of the existence part of Theorem \ref{existence and uniqueness}.
Now the result follows from the uniqueness
part of Theorem \ref{existence and uniqueness}.
\hfill $\square$

\stm\label{Example of big open not everything}
{\bf Example where a canonical reduction is possible only on $U \ne X$.}
This example is well known to experts (we learnt it from 
Tomas Gomez), but not easily accessible in the literature. 
Let $k$ be any field, $X = \PP^2_k = \proj k[x,y,z]$ with the usual $\OO_X(1)$,
and $G = GL_2$ together with the diagonal subgroup as the chosen maximal torus
and the upper triangular subgroup as the chosen Borel. Let $P = (0,0,1)$, 
which is a closed point in $X$ with residue field $k$, and let
$\I_P \subset \OO_X$ be its ideal sheaf. We now make a series of claims, each
of which is essentially an exercise which we leave to the reader.

(1) $Ext^1(\I_P,\OO_X(1)) \cong k$, hence there exists a 
non-split short exact sequence
$0 \to \OO_X(1) \stackrel{i}{\to} \E \stackrel{j}{\to}\I_P \to 0$
of $\OO_X$-modules, and the isomorphism class of
an $\OO_X$-module $\E$ that occurs in such a non-split exact sequence is unique.
Such an $\E$ is a locally free $\OO_X$-module. 

(2) The injective homomorphism of sheaves 
$i : \OO_X(1) \to \E$ is fiberwise injective on the big open subset $X - \{ P\}$
of $X$, while it is $0$ on the fiber at $P$.

(3) The filtration $0\to \OO_X(1) \stackrel{i}{\to} \E$ is the
Harder-Narasimhan filtration
of $\E$ both in the $\mu$-semistability sense (which concerns us here) and the
Gieseker-semistability sense. 

(4) Let $E$ be the uniquely determined principal $GL_2$-bundle corresponding to
$\E$. Then the above filtration gives the canonical reduction of $E$, and it
is defined on the open subset $X - \{ P\}$ but cannot be extended to all of $X$.

\section{Families and Relative reductions}

As before, $G$ is a split reductive group over a field $k$ with a 
chosen split maximal torus and a Borel containing it, which satisfies
the hypothesis $(*)$.

\stm\label{Families of varieties and bundles on them}
{\bf Families of varieties and bundles on them.}
From now onwards, instead of an extension field $K$ of $k$ and
a geometrically irreducible smooth
projective variety $X$
over $K$ equipped with a very ample line bundle $\OO_X(1)$,
we consider a noetherian 
scheme $S$ over $k$ and a smooth proper morphism $\pi: X\to S$
with geometrically connected
fibers of a constant dimension $n$, equipped with a relatively very ample line bundle $\OO_{X/S}(1)$.
An open subscheme $U\subset X$ is said to be {\bf relatively big} over $S$ 
if the fiber $U_s$ is big in $X_s$ for each $s\in S$.
Instead of just
considering $G$-bundles on the projective varieties
$X_L$ for field extensions $L/K$, we will consider
arbitrary (not necessarily noetherian)
$S$-schemes $T$ and $G$-bundles on $X_T = X\times_S T$.
We denote by $\OO_{X_T/T}(1)$ the pullback of $\OO_{X/S}(1)$ to $X_T$.
The category of all $T$-schemes
will be equipped with the fppf topology.

\begin{lemma}\label{Finiteness of the set of types over S} 
{\bf Finiteness of the set of types over $S$}
Let $\pi : X\to S$ and $\OO_{X/S}(1)$ be as in
(\ref{Families of varieties and bundles on them}) above where
$S$ is a noetherian $k$-scheme,
and let $E$ be a principal $G$-bundle on $X$.
For any $s\in S$, let $\HN(E_s) \in \ov{C}$ denote the canonical type
of the restriction $E_s = E|_{X_s}$ where $X_s$ is the
schematic fiber of $\pi$ over $s$, equipped with the very ample line bundle
$\OO_{X_s}(1) = \OO_{X/S}(1)|_{X_s}$. Then the set
$\{ \HN(E_s)\,|\, s\in S \}$
is finite.
\end{lemma}

{\bf Proof.} Note that the set of all standard 
parabolics $P$ in $G$ is finite, and for any such $P$,
there exists a finite set of dominant virtual characters
(the basic positive virtual characters of
Lemma \ref{characters positive combination of simple roots}) 
$\chi_{I_P,\alpha}$ on $P$ where $\alpha$ varies over $\Delta - I_P$,
such that any dominant characters $P \to \GG_m$ is a non-negative
rational linear combination of the $\chi_{I_P,\alpha}$. We can choose
a common integer $N>0$ such that each $N\chi_{I_P,\alpha} \in X^*(T)$.
For simplicity, we will use the notation $\ov{\alpha} =  N\chi_{I_P,\alpha}$.

As $S$ is noetherian, it is a finite disjoint union of connected open
subschemes, so we can without loss of generality assume that $S$ is connected.
Let $G \to GL(V_{\ov{\alpha}})$ denote the representation defined by the 
dominant weight  $\ov{\alpha} = N\chi_{I_P,\alpha}$.
Let $E(V_{\ov{\alpha}})$ denote the vector
bundle on $X$ associated to the principal $G$-bundle $E$ via the representation
$V_{\ov{\alpha}}$.

As $S$ is noetherian, the set of Harder-Narasimhan types (in the sense of
Gieseker semistability) of the
vector bundles $E_s(V_{\ov{\alpha}})$, as $s$ varies over $S$, is a finite set.
Hence as $s$ varies over $S$,
the set of degrees of all line bundles $L$ on $X_s$ for which there exists at least  
one injective homomorphism of sheaves of
$\OO_{X_s}$-modules $i : L \to E_s(V_{\ov{\alpha}})$ is bounded 
above by a constant $c_{\ov{\alpha}}$.
If $(P,\sigma)$ is a reduction of $E_s$, then applying this
to the resulting homomorphisms
$i: L_{\sigma, \ov{\alpha}} \to E_s(V_{\chi_{\ov{\alpha}}})$
where $\alpha \in \Delta - I_P$, we see that
$$\langle N\chi_{I_P, \alpha}, \mu_{(P,\sigma)}(E_s) \rangle
= \deg(L_{\sigma, \ov{\alpha}}) \le c$$
where $c$ is the maximum of the $c_{\ov{\alpha}}$ as $\alpha$ varies over
$\Delta - I_P$. 
If $\beta \in I_P$, then $\langle \beta, \mu_{(P,\sigma)}(E_s) \rangle =0$ while
for any $\theta \in \wh{G}$, we have 
$$\langle \theta, \mu_{(P,\sigma)}(E_s)\rangle =\langle \theta, \mu_{(G,\id)}(E_s)\rangle,$$
and $\mu_{(G,\id)}(E_s)$ is independent of $s$ as $S$ is connected.
The set $A$ of types of the $E_s$ as $s$ varies over $S$ 
is contained in the intersection of
the closed positive Weyl chamber $\ov{C}$ with
a lattice $(1/n_0) X_*(T)$ in $\Q\otimes X_*(T)$ 
by Lemma \ref{types lattice} where the integer $n_0$ depends only on $G$,
and is independent of $s$.
Hence the desired conclusion follows by Lemma \ref{boundedness and finiteness}.
\hfill$\square$

\begin{lemma}
Let $X\to S$ and $\OO_{X/S}(1)$ be as above, with $S$ noetherian, and let
$E$ be a principal $G$-bundle on $X$. Suppose that there exists a 
open subscheme $U\subset X$ that is relatively big over $S$, and 
a section $\sigma : U \to E/P$ where $P$ is a standard parabolic in $G$.
For any $s\in S$, 
let $U_s = U\cap X_s$ and let $\sigma_s= \sigma|_{U_s}: U_s \to E_s/P$.
Then $S$ is the disjoint union of finitely many open subschemes
$S_i$ such that the type of $(P,\sigma_s)$ is constant over $s\in S_i$.
\end{lemma}

{\bf Proof.} Note that if $\L$ is a line bundle on $U$, then the function
$s\mapsto \deg(\L|_{U_s})$ is locally constant on $S$ as $U$ is relatively big over $S$.
For each $\ov{\alpha} = N\chi_{I_P,\alpha}: P \to \GG_m$ as in the
notation in the proof of the Lemma \ref{Finiteness of the set of types over S},
let $\L_{\ov{\alpha}}$ be the line bundle
$\ov{\alpha}_*\sigma^*E$ on $U$. As $U$ is relatively big over $S$,
$\deg(\L|_{U_s})$ is locally constant over $S$. As $\alpha$ varies over the
finite set $\Delta - I_P$, we have a decomposition of $S$ as a
finite disjoint union of open subschemes $S_i$ such that the function
$s\mapsto \deg(\L|_{U_s})$ is constant on $S_i$ for all $\alpha$.
The result follows by arguing as in the proof of
the Lemma \ref{Finiteness of the set of types over S}.
\hfill$\square$

\stm{\bf Definition of a relative reduction.}
With notation as in (\ref{Families of varieties and bundles on them}) above,
we recall from [G-N-2]
the definition of a relative reduction of a $G$-bundle
$E$ on $X_T$ to a standard parabolic $P\subset G$.
An open subscheme $U\subset X_T$ is called
{\it relatively big} over $T$ if $U_t$ is big in
$X_t$ for all $t\in T$.
Let $\sigma : U \to E/P$ be a section of $E/P \to T$.
We say that two such pairs $(U_1,\sigma_1)$ and $(U_2,\sigma_2)$ are equivalent
if $\sigma_1|_{U_1\cap U_2} = \sigma_2|_{U_1\cap U_2}$, which can be
seen to be an equivalence relation. We denote the 
equivalence class of $(U,\sigma)$ by $[U,\sigma]$.

Let $E(V_P)$ denote the vector bundle on $X$  
associated to the principal $G$-bundle $E\to X$ via  
the linear $G$-representation $V_P = H^0(G/P, \omega^{-1}_{G/P})^{\vee}$.
Consider its restriction $E(V_P)|_U$ to the relatively big open subscheme $U$.
Recall from (\ref{definition of V_P})
the $P$-invariant $1$-dimensional linear subspace $J\subset V_P$, 
on which $P$ acts by the dominant character
$\lambda_P: P\to \GG_m$.
As in (\ref{Alternative description of reduction}), 
the reduction $\sigma : U \to E/P$ gives rise to a line subbundle
$L_U \subset E(V_P)|_U$. The line bundle $L_U$ equals $(\lambda_P)_* \sigma^* E$.
We say that $(U,\sigma)$ defines a relative reduction
of the structure group of $E$ from $G$ to $P$ is there exists a
prolongation $(L,\phi)$ of $L_U$ to $X$, where
$L$ is a line bundle on $X$ and $\phi : L|_U\to L_U$ is an isomorphism. 

If such a prolongation $(L,\phi)$ exists then is unique up to unique isomorphism,
and comes with a uniquely defined
injective homomorphism of $\OO_X$-modules $i: L \to E(V_P)$ which prolongs
the inclusion $L_U \hra E(V_P)|_U$. For a given $(L,\phi)$, the pair $(L,i)$ is 
unique up to unique isomorphism. This is a consequence of the relative bigness
of $U$ in $X/S$, as shown in [G-N-2] (see the Introduction and the
Proposition 2.2 and Lemma 2.3 of [G-N-2]). Here, we say that a pair $(L,i)$ is
isomorphic to another pair $(L',i')$ if there exists an isomorphism
$\phi : L\to L'$ such that $i = i'\circ \phi$.
The isomorphism class $[L,i]$ 
of $(L,i)$ is uniquely determined by the equivalence class $[U,\sigma]$ of $(U,\sigma)$.

Conversely, 
consider all pairs $(L,i)$ where $L$ is a line bundle on $X$ 
and $i: L\to E(V_P)$ is an injective $\OO_X$-linear 
homomorphism of sheaves,
such that \\
(i) the open subscheme $U = \{ x\in X\,|\,\rank(i_x) =1\} \subset X$ 
is relatively big over $S$, that is, for each $s\in S$
the fiber $U_s$ has complementary codimension 
$\ge 2$ in the fiber $X_s$, and \\
(ii) the section $U\to {\bf P}(E(V_P))$  
defined by $i$ factors via the closed embedding 
$E/P \hra {\bf P}(E(V_P))$ given by the very ample line bundle
$\omega^{-1}_{G/P}$, to define a section $\sigma: U\to E/P$.\\ 
An isomorphism of pairs $(L,i) \to (L',i')$
is an isomorphism of line bundles 
$\phi :L \to L'$ such that $i = i'\circ \phi$. 
Note that there exists at most one isomorphism $\phi$ between any two pairs,
so the groupoid of all such pairs is equivalent to the set of all their
isomorphism classes.
As in the absolute case (see statement 
\ref{recovering the reduction from the homomorphism i} above)  
the pair $(L,i)$ arises from the pair $(U,\sigma)$ where $U$ and $\sigma$
are as defined in (i) and (ii) above. This correspondence establishes
a bijection between (a) the set of all equivalence classes $[U,\sigma]$ of 
pairs $(U,\sigma)$ for which $L_U$ admits a prolongation to $X$, and (b)
the set of all isomorphism classes $[L,i]$ of pairs $(L,i)$ which satisfy
(i) and (ii) above.

Finally, as in [G-N-2], we define a rational reduction of the structure group of $E$ to
a parabolic $P\subset E$ to be an isomorphism class of pairs $(L,i)$
that satisfy (i) and (ii) above. Equivalently, such a reduction is an equivalence
class of pairs $(U,\sigma)$ where $U\subset X$ is relatively big over $S$ and 
$\sigma : U\to E/P$ is a section of $E/P \to X$ defined over $U$, such that the
resulting line bundle $L_U$ on $U$ admits at least one prolongation to $X$.

In the special case where $T = \spec K$ for an extension field $K/k$,
the above definition is equivalent to the usual definition
of a rational reduction to $P$ (as defined by $\sigma :U \to E/P$)
that we recalled earlier
(see Proposition 3.4 in [G-N-2]).

\stm\label{Local constancy of type}{\bf Local constancy of type.}
If $E$ is a bundle on $X_T$ and if
there exists a relative reduction $(U,\sigma)$ of $E$ to a parabolic $P$,
then note that the set-valued function $t\mapsto \mu_{(P,\sigma_t)}(E_t)$ is locally 
constant on $t$. 

\stm{\bf Relative canonical reduction of constant type $\tau$.}
Finally, for any $\tau \in \ov{C}$, we say that a pair 
$(L,i)$ as above defines a 
canonical reduction $[L,i]$ of constant type $\tau$ if its restriction to each
fiber $X_t$ of $X_T\to T$ is a canonical reduction of $E_t$ 
of constant type $\tau$.
By Theorem \ref{existence and uniqueness}, the 
type of a bundle
is well-defined under the hypothesis ${\bf (*)}$ of the Introduction.

\stm{\bf Descent for relative reductions.}
If $S' \to S$ is an fppf morphism and if we have two reductions $(L_1,i_1)$
and $(L_2,i_2)$ of $E$ over $S$ whose pullbacks to $S'$ are isomorphic
then the two reductions are isomorphic on $S$ itself, and the isomorphism
on $S$ (which is unique when it exists) pulls back to the isomorphism
over $S'$ (which is unique when it exists). This means fppf descent
holds for relative reductions.

The following lemma shows the effectivity of fppf descent, that is, 
if $S' \to S$ is an fppf morphism and we have a relative rational 
$P$-reduction $[L',i']$ 
of $E_{S'}/X_{S'}/S'$ such that under the two projections
$p_1,p_2: S'\times _S S' \toto S$ the pull backs $p_1^*(L', i')$
and $p_2^*(L', i')$ are equivalent reductions,
then there exists a unique relative rational reduction $[L,i]$ of $E/X/S$
to $P$ which pulls back to $[L',i']$.

\begin{lemma}\label{sheaf property}{\bf (Effective fppf descent for reductions.)}
Let $Y$ be a scheme, $\E$ be a sheaf of $\OO_Y$-modules, 
and $p: Y'\to Y$ be an fppf cover of $Y$. Let $L'$ be a line bundle on $Y'$
and let there be given an injective $\OO_{Y'}$-linear homomorphism 
$i': L' \to p^*\E$ such that there exists $g \in \GG_m(Y'\times_YY')$ such that 
$p_1^*(i') = g p_2^*(i')$ (we do not assume a $1$-cocycle condition on $g$). 
Then there exists a line bundle $L$ on $Y$ and an injective 
$\OO_Y$-linear homomorphism $i: L \to \E$, such that
$(L',i' : L' \to p^*\E)$ is the pullback of $(L, i)$.
The pair $(L,i)$ is unique up to a unique isomorphism.
\end{lemma}

\proof The image subsheaf $\im(i')
\subset \E'$ is locally free of rank $1$ as $i'$ is injective. By the existence
of $g\in \GG_m(Y'\times_YY')$, it 
satisfies $p_1^*(\im(i')) = p_2^*(\im(i'))$, and so being a subsheaf
of the pullback $\E$, the $1$-cocycle condition is automatically satisfied.
Hence by effective fppf descent for quasi-coherent sheaves,
we obtain a line bundle $L$ on $Y$ with an
injective homomorphism $i: L\to \E$ together with an isomorphism $p^*(L,i) = (L',i')$. 
By descent, the pair $(L,i)$ is unique up to a unique isomorphism.
\hfill $\square$

{}

\section{The algebraic stacks of canonical types}

Let $\pi: X\to S$ be a proper morphism of noetherian schemes, 
and let $\F$ be a coherent $\OO_X$-module that is flat over $S$.
For any $S$-scheme $T$ and a point $t\in T$, we have a $\kappa(t)$-linear
base-change
homomorphism $({\pi_T}_* \F_T)\otimes_{\OO_T}\kappa(t) \to
H^0(X_t,\F_t)$ where $\kappa(t)$ is the residue field at $t$.
Hence given any invertible sheaf $M$ of $\OO_T$-modules and any
$\OO_T$-linear homomorphism $j : M \to  {\pi_T}_*\F_T$,
we define the $\kappa(t)$-linear
{\it evaluation homomorphism} $j_t: M \otimes_{\OO_T}\kappa(t) \to H^0(X_t,\F_t)$
to be the composite
$$M \otimes_{\OO_T}\kappa(t) \stackrel{j\otimes \id}{\to}
({\pi_T}_*\F_T) \otimes_{\OO_T}\kappa(t)\to H^0(X_t,\F_t).$$

{\bf The functor $\Psi_{\F/X/S}$.} For any $S$-scheme $T$, consider pairs
$(M, i)$ where $M$ is an invertible sheaf of $\OO_T$-modules 
and $i : M \to {\pi_T}_*\F_T$ is an $\OO_{X_T}$-linear homomorphism such that
at each $t\in T$, the evaluation 
homomorphism $i_t: M \otimes_{\OO_T}\kappa(t) \to H^0(X_t,\F_t)$ is non-zero.
We say that two such pairs $(M,i)$ and $(M',i')$ are equivalent if
there exists an $\OO_T$-linear isomorphism $g: M \to M'$ such
that $i = i'\circ g$. This is clearly an equivalence relation, and all equivalence classes
form a set $\Psi(T)$. For any morphism $f: T' \to T$, by defining
$f^*(M,i)$ to be $(f^*M, f^*i)$ we get a well defined map of sets
$f^*: \Psi(T) \to \Psi(T')$. This defines a contravariant functor
$\Psi$ from $S$-schemes to sets.
We will denote this functor by $\Psi_{\F/X/S}$ for clarity when required.

The following lemma is a projective 
version of the result of Grothendieck on the 
representability by a linear scheme for sections of direct images 
(see [EGA III 7.7.8, 7.7.9], and [N-2, 5.8] for an exposition). 
We expect this lemma, though elementary, to be of independent interest.

\begin{lemma}\label{projective representation}
Let $X\to S$ be a proper morphism of noetherian schemes, and let
$\F$ be a sheaf of $\OO_X$-modules that is flat over $S$. Then the following holds.

(1) There is a natural isomorphism of functors between the functor of points
of $Proj_S \, Sym_S ^{\bullet}(\QQ)$,
where $\QQ$ denotes the Grothendieck Q-sheaf of $\F/X/S$, and the 
above functor $\Psi_{\F/X/S} : (Schemes/S)^{op} \to Sets$, defined by
applying $\un{Hom}(-,\OO_T)$ over any $S$-scheme $T$.

(2) Consequently, the functor $\Psi_{\F/X/S}$
is representable by the $S$-scheme 
$\PP(\QQ ) = Proj_S \, Sym_S ^{\bullet}(\QQ)$, and 
a universal element (Poincar\'e family) for $\Psi_{\F/X/S}$
parameterized by $\PP(\QQ )$ is of the form
$(\OO_{\PP(\QQ)}(-1), i)$,
where the homomorphism $i: \OO_{\PP(\QQ)}(-1) \to (\pi_{\PP(\QQ)})_* \F_{\PP(\QQ)}$
is obtained by applying $\un{Hom}(-, \OO_{\PP(\QQ)})$ to the
universal quotient 
$q: \pi_{\PP(\QQ)}^*\QQ \to \OO_{\PP(\QQ)}(1)$.
\end{lemma}

{\it Proof.} If $\E$ is any coherent sheaf on $S$, 
then recall that 
$\PP(\E ) = Proj_{S}Sym_S^{\bullet}(\E)$ represents the functor 
$\varphi : (Schemes/S)^{op} \to Sets$ which associates to 
to any $S$-scheme $T$ the set of all equivalence classes of pairs $(L,q)$ 
where $L$ is an invertible sheaf on $T$ and $q: \E_T \to L$
is a surjective $\OO_T$-linear homomorphism. Two such pairs $(L,q)$ and $(L',q')$
are defined to be equivalent if there exists an $\OO_T$-linear isomorphism
$\lambda : L\to L'$ such that $\lambda \circ q = q'$.

We will apply the above by taking the coherent $\OO_S$-module $\E$ to be
the Grothendieck sheaf $\QQ$ of the sheaf $\F$ on $X$ which is flat over $S$.
Recall that locally over $S$, 
we can take $\QQ$ to be the cokernel of the transpose of the $0$ th 
differential of a Grothendieck semicontinuity complex for $\F/X/S$.
Though the semicontinuity complexes are defined only locally over $S$, these cokernels
uniquely glue to define $\QQ$ globally on $S$. 
Note that as a consequence of the base change property for a semicontinuity complex,
the Q-sheaf has the base change property, that is, if $f: T \to S$ is any
morphism, then the Q-sheaf $\QQ_T$ for
$\F_T/X_T/T$ is the pull back $\QQ_T = f^*\QQ$ under $f$ of the Q-sheaf $\QQ$
for $\F/X/S$.

Recall that the Grothendieck sheaf $\QQ$ of $\F$ has the universal property 
(see [EGA III 7.7.8, 7.7.9]) that for any $S$-scheme $T$ we have a natural bijection
$$Hom(\QQ_T,\OO_T) \stackrel{\sim}{\to}\Gamma(T,{\pi_T}_*\F_T)$$
which is in fact an isomorphism of $\Gamma(T,\OO_T)$-modules.
As this isomorphism 
is functorial in $T$, 
by pull-back to any point $t$ of $T$ it 
follows that surjective homomorphisms $\QQ_T \to \OO_T$ exactly 
correspond to sections $s\in \Gamma(T,{\pi_T}_*\F_T)$
for which the corresponding sections $s_t \in H^0(X_t, \F_t)$
obtained by base change are non-zero for each $t\in T$.

Let $T$ be any $S$-scheme, and let 
$h : \QQ_T \to L$ be a surjection where $L$ is an invertible $\OO_T$-module.
Let $U_i$ be an open cover of $T$ over which $L$ has a local basis $e_i$
with transition functions $g_{ij} \in \GG_m(U_i\cap U_j)$ with $e_ig_{ij} = e_j$.
Then $h$ corresponds to a family of homomorphisms $h_i : \QQ_{U_i} \to \OO_{U_i}$
with $h_i = g_{ij}h_j$, such that each $h_i$ is surjective at all points of $U_i$. 
Hence $h_i$ corresponds to a section $s_i \in \Gamma(U_i, (\pi_{U_i})_*\F_{U_i})$
with $s_i = g_{ij}s_j$ such that each $s_i$
evaluates under base change to give a non-zero element
$(s_i)_t \in H^0(X_t,\F_t)$ at each $t\in U_i$.

The dual line bundle $L^*$ has local bases $\sigma_i$ over $U_i$ with
$\langle e_i, \sigma_i\rangle = 1$, and so $\sigma_i = g_{ij}\sigma_j$.
Hence the $s_i$'s glue together to
define an $\OO_T$-linear homomorphism $s: L^* \to {\pi_T}_*\F_T$,
with $\sigma_i \mapsto s_i$, which
evaluates at each $t\in T$ to give a non-zero $\kappa(t)$-linear
map $L^*\otimes_{\OO_T} \kappa(t) \to H^0(X_t, \F_t)$. Hence
$s \in \Psi(T)$. Note that $s$ is obtained by applying
$\un{Hom}(-,\OO_T)$ to $h: \QQ_T \to L$.

Conversely, given any
$s: L^* \to {\pi_T}_*\F_T$ in $\Psi(T)$, we can re-trace the above steps to obtain
a surjection $h : \QQ_T \to L$. This establishes an isomorphism of
functors between the functor represented by $\PP(\QQ)$ and the functor $\Psi$, 
defined by applying $\un{Hom}(-,\OO_T)$ over any $S$-scheme $T$. 
Applying this isomorphism to the universal family (tautological quotient line bundle)
$q: \pi_{\PP(\QQ)}^*\QQ \to \OO_{\PP(\QQ)}(1)$ on $\PP(\QQ)$ gives the
desired universal family for $\Psi$ parameterized by $\PP(\QQ)$.
\hfill$\square$

\rem
Let $\pi: X\to S$ be a projective morphism, with $S$ noetherian.
Unlike the functor $T \mapsto \Gamma(T,{\pi_T}_*\F_T)$, 
which is representable if and only if $\F$ is flat over $S$ (see [N-1]),
the functor $\Psi$ may or may not be representable when $\F$ 
is not flat. For example, let $X = S = \spec k[t]$ for a field $k$, and let
$\F_1 = (k[t]/(t))^{\sim}$ and $\F_2 = \F_1 \oplus \OO_X$,   
which are coherent sheaves on $X =S$ which are not flat over $S$.
We now show that 
the corresponding functor $\Psi_{\F_1/X/S}$ is representable, 
while $\Psi_{\F_2/X/S}$ is not representable.
For this, we can assume $k= \ov{k}$.
For $\F_1$, the scheme $\spec k[t]/(t)
\hra \spec k[t]$ represents $\Psi$. For $\F_2$, suppose $\Psi$ is 
represented by a scheme $p: R\to \spec k[t]$. The fiber of 
$R$ over the origin 
$x_0 = \spec k[t]/(t)$ must be the projective line $p^{-1}(x_0) = \PP^1_k$,
while over $U = \spec k[t,t^{-1}]$ the restriction $R_U \to U$ is
an isomorphism. Let $V \subset R$ be an affine open
subscheme with $V \cap p^{-1}(x_0) \ne \emptyset$. Then as
$k= \ov{k}$, we will have at least two distinct $k$-valued points $y_1$ 
and $y_2$ in $V \cap p^{-1}(x_0)$. From the definition of 
$\F_2$, it can be seen that there will exist global sections
$s_1,s_2 \in \Gamma(S,R)$ with $s_i(x_0) = y_i$. But then 
$s_1|_U = s_2|_U$, which contradicts the separatedness of  
$V \to S$.

\begin{lemma}\label{Blaah}
Let $Y$ be a smooth irreducible projective variety over a field $K$
with a very ample line bundle $\OO_Y(1)$.
Let $\E$ be a vector bundle on $Y$ and let $L$ and $M$ be line bundles on
$Y$.  Let $i : L \to \E$ and 
$ h : M \to \E$ be non-zero $\OO_Y$-homomorphisms such that\\
(i) inside the fiber of $\E$ at the generic point of $Y$, the image of
$i$ and the image of $h$ are the same $1$-dimensional linear subspace, and\\
(ii) $h$ is fiberwise injective on a big open subscheme $U$ of $Y$.

Then we must have $\deg(L) \le \deg(M)$.
Moreover, if $\deg(L) = \deg(M)$ then there exists an isomorphism
$\phi : L \to M$ with $i = h \circ \phi$.
\end{lemma}


\proof Recall that the {\it saturation} $\ov{F}$ of a coherent subsheaf $F$
of $\E$ is the inverse
  image $q^{-1}(T)$ under the quotient homomorphism $q: \E \to \E/F$ of the torsion
  subsheaf $T$ of $\E/F$. Let $\ov{L} \subset \E$ be the saturation of $\im(i)$ in $\E$.
The condition (ii) implies that $\ov{M} = \im(h)$ is already saturated in $\E$
  (means its saturation is itself). 
By (i), we have $\ov{L} = \ov{M} \subset \E$. As $h: M \to \ov{M}$ is an isomorphism, 
  we have $\deg(M) = \deg(\ov{M})$. As $L$ is a subsheaf
  of $\ov{L}$ generically equal to it, we have $\deg(L) \le \deg(\ov{L})$. 
  Hence we get the inequality  $\deg(L)\le  \deg(M)$.
  The equality $\deg(L) = \deg(\ov{L})$ holds if and only if the set of points
  where $i : L\to \E$ is fiberwise injective form a big open subscheme $W$. 
  Then the intersection $U\cap W$, where both $h$ and $i$ are fiberwise injective,
  is a big open in $Y$, and $L|_{U\cap W} = M|_{U\cap W} \subset \E|_{U\cap W}$.
  If $\phi' : L|_{U\cap W} \to M|_{U\cap W}$ is the induced isomorphism, then as
  $U\cap W$ is big in $Y$, $\phi'$ has a unique prolongation to an isomorphism
  $\phi : L\to M$ which has the desired property. \hfill$\square$

\begin{theorem}\label{valuative criterion}
{\bf Valuative criterion for canonical type.}
Let $S = \spec R$ where $R$ is a discrete valuation ring over $k$, with generic
point $s_0$ and closed point $s_1$.
Let $\pi : X\to S$ be a smooth projective morphism with
geometrically irreducible fibers, and $\OO_{X/S}(1)$ a relatively very ample
line bundle on $X$ over $S$, and let 
$X_i = \pi^{-1}(s_i)$ for $i = 0,1$. Let $E$ be a principal $G$-bundle on $X$.
Then we have the following.\\
(a) $\HN(E_0) \le \HN(E_1)$, \\
(b) $\HN(E_0) = \HN(E_1)$ only if a relative canonical reduction for
$E$ exists over $S$.\\
In particular if $E_1$ is semistable then so is $E_0$, and there exists
a relative canonical reduction for $E$ over $S$.
\end{theorem}  

\proof The proof will be divided into a sequence of steps.

{\it Step (1)}: Let $\sigma_0 : U_0 \to E_0/P$ be the canonical reduction of
$E_0$, where $U_0$ is the maximal big open subscheme of $X_0$
for the reduction. Then there exists 
an open subscheme $U\subset X$ and a reduction
$\sigma : U \to E_U/P$ such that $U_0 = U\cap X_0$, $\sigma_0 = \sigma|_{U_0}$, 
and $U_1 = U\cap X_1$ is nonempty.

{\it Proof of Step (1)}: Recall that if $\pi: Y\to Z$ is a proper morphism of schemes
where $Z$ is noetherian, integral and regular in codimension $1$, then any generic section
of $\pi$ admits a unique maximal prolongation to a section defined over an open subscheme
$U$ of $X$, which is moreover a big open subscheme of $Z$.
We will apply this to the projection $E/P \to X$. 
Note that $X$ is a separated noetherian regular
integral scheme, $X_0$ is a dense open subscheme of
$X$ and $X_1$ is an effective divisor in $X$. 
In particular, any big open subset of $X$ has a nonempty intersection with $X_1$.
The
section $\sigma_0 : U_0 \to E/P$ has a unique largest 
extension to a section $\sigma: U \to E/P$ defined over a big open subscheme
$U\subset X$. Therefore, $U_0 = U\cap X_0$ with $\sigma_0 = \sigma|_{U_0}$,
and the intersection $U_1= X_1\cap U$ is nonempty.

{\it Step (2)}: Let $\sigma_1 : U_1 \to E_1/P$ be the restriction of $\sigma$ to $U_1$,
which is a rational reduction of $E_1$ to $P$ as $U_1\ne \emptyset$. 
It has a unique maximal extension
$\sigma_1' : U_1'\to E_1/P$ with $U\subset U_1'\subset X_1$
where $U_1'$ is a big open subscheme of $X_1$.
Then
$$\mu_{(P,\sigma_1')}(E_1) \ge \mu_{(P,\sigma_0)}(E_0)$$

{\it Proof of Step (2).} To prove the above inequality, we have to equivalently
prove that for all dominant characters $\chi : P\to \GG_m$, we must have
$\langle \chi, \mu_{(P,\sigma_1')}(E_1)\rangle
\ge \langle \chi, \mu_{(P,\sigma_0)}(E_0)\rangle$.
By definition of the type of a reduction,
this is the inequality
$$\deg_{X_1}\chi_*(\sigma_1')^*E_1 \ge  \deg_{X_0}\chi_*(\sigma_0)^*E_0.$$
Let $\chi_*\sigma^*E$ be the line bundle on $U$ associated
to the reduction $\sigma : U \to E/P$. As $U$ is a big open in $X$ and $X$ is
a separated noetherian regular integral scheme, this line bundle
admits a unique prolongation $L_{\sigma, \chi}$ to $X$, together with 
a homomorphism $i : L_{\sigma, \chi} \to E(V_{\chi})$
as in (\ref{Alternative description of reduction}), which is fiberwise injective
on $U$. For $j=0,1$, let $L_j = L|_{X_j}$ and $i_j = i|_{X_j} : L_j \to E_j(V_{\chi})$ 
which is therefore fiberwise injective on $U_j$. 

There exists a maximal big open subset $U_1'\subset X_1$ on which
the generic section $\sigma_1 : U_1\to E_1/P$ prolongs, defining a section
$\sigma_1' : U_1' \to E_1/P$ with $\sigma_1'|_{U_1} = \sigma_1$. 
It defines a line bundle $\chi_*(\sigma_1')^*E_1$ on $U_1'$, which
as in (\ref{Alternative description of reduction}) admits a unique
prolongation to a line bundle $M_1$ on $X_1$ together with a 
homomorphism $f : M_1 \to E_1(V_{\chi})$ that is fiberwise injective on $U_1'$.
As $\sigma_1'|_{U_1} = \sigma_1$, we have $M_1|_{U_1} = L_1|_{U_1}$,
and $f|_{U_1} = i_1|_{U_1}$. 
Hence it follows by Lemma \ref{Blaah} that
$\deg_{X_1}(M_1) \ge \deg_{X_1}(L_1)$.
Now,
$$\deg_{X_1}(M_1) = \deg_{X_1}\chi_*(\sigma_1')^*(E_1) = 
\langle \chi, \mu_{(P,\sigma_1')}(E_1)\rangle$$
while
$$\deg_{X_1}(L_1) = \deg_{X_0}(L_0) =
\deg_{X_0}\chi_*(\sigma_0)^*(E_0) =
\langle \chi, \mu_{(P,\sigma_0)}(E_0)\rangle .$$
Hence $\langle \chi, \mu_{(P,\sigma_1')}(E_1)\rangle 
\ge \langle \chi, \mu_{(P,\sigma_0)}(E_0)\rangle$ as desired, proving Step (2).

{\it Step (3)}: If $\mu_{(P,\sigma_1')}(E_1) = \mu_{(P,\sigma_0)}(E_0)$
then there exists a relative reduction of $E$ to $P$ on $X/S$, which
restricts to $(P, \sigma_0)$ over $s_0$ and restricts to $(P, \sigma_1')$ over $s_1$.

{\it Proof of Step (3)}: If $\mu_{(P,\sigma_1')}(E_1) = \mu_{(P,\sigma_0)}(E_0)$
then for each dominant $\chi: P\to \GG_m$, with notation as in the
proof of Step (2), we have 
$\deg_{X_1}(M_1) =
\langle \chi, \mu_{(P,\sigma_1')}(E_1) \rangle =
\langle \chi, \mu_{(P,\sigma_0)}(E_0) \rangle = \deg_{X_0}(L_0)
= \deg_{X_1}(L_1)$. Hence by Lemma \ref{Blaah},
we have an isomorphism $\phi : L_1 \to M_1$ with $i_1 = h\circ \phi$.
Hence $U_1 = U_1'$, which is  big in $X_1$,
showing that $U\subset X$ is relatively
big over $S$. The line bundle 
$(\lambda_P)_*\sigma^*E$ on the big open subscheme $U$ of $X$ prolongs 
to a line bundle on $X$ as $X$ is a
noetherian integral separated regular scheme.
Hence $\sigma : U \to E/P$
is a relative reduction of $E$ over $S$, proving Step (3). 

With the above steps, we are now ready to prove the statements (a) and (b)
of the theorem.

{\it Proof of (a)}: By Step (2), $\mu_{(P,\sigma_1')}(E_1) \ge \mu_{(P,\sigma_0)}(E_0)$.
As by Corollary \ref{maximality of HN type}
$\HN(E_1)$ is the unique maximum in $\Q\otimes X_*(T)$ for types of reductions of $E_1$ 
to standard parabolics, we have $\HN(E_1) \ge \mu_{(P,\sigma_1')}(E_1)$.
As $\HN(E_0) = \mu_{(P,\sigma_0)}(E_0)$, it follows that
$\HN(E_1) \ge \HN(E_0)$ as claimed.

{\it Proof of (b)}: If $\HN(E_1) = \HN(E_0)$, then as
$\HN(E_1) \ge \mu_{(P,\sigma_1')}(E_1) \ge \mu_{(P,\sigma_0)}(E_0)
= \HN(E_0)$, we have $\mu_{(P,\sigma_1')}(E_1) = \mu_{(P,\sigma_0)}(E_0)$. 
Hence by Step (3), we have a relative reduction $\sigma : U\to X$ of $E$ over $S$,
and it restricts to the canonical reductions on both $X_1$ and $X_0$.
This completes the proof of (b) and of the theorem.
\hfill$\square$

\begin{lemma}\label{exercise} Let $S$ be a noetherian scheme and let $C\subset |S|$
be a constructible subset. Suppose $C$ satisfies the property that if
$R$ is any discrete valuation ring and $f: \spec R\to S$ a morphism
such that $f$ maps the generic point of $\spec R$ to a point of $C$, then $f$
maps the closed point of $\spec R$ to a point of $C$. Then $C$ is a closed subset
of $|S|$.
\end{lemma}

{}

\proof Any constructible subset $C$ contains a dense open subset of its closure in $S$. Hence
to show that $C$ is closed, it is enough to show that any 
specialization $x\in S$ of any point $y\in C$ is again in $C$. Let $Z$ be the closure
of $y$ in $S$ with reduced subscheme structure. The local ring $\OO_{Z,x}$ is a noetherian
local domain, hence (see e.g. [Stacks Project] Lemma 10.118.13)
there exists a dvr $R$ and a finite surjective morphism
$\spec R \to \spec \OO_{Z,x}$ that sends the generic point to the generic point and
the closed point to the closed point. By applying the hypothesis
of the lemma to the composite $\spec R \to \spec \OO_{Z,x} \hra S$, it follows that $x\in C$.
\hfill $\square$

\begin{proposition}\label{openness of semistability}
{\bf (Openness of semistability.) }
Let $X\to S$ and $\OO_{X/S}(1)$ be as above with $S$ noetherian, 
and let $U\subset X$ be an open subscheme
that is relatively big over $S$. Let $E$ be a principal $G$-bundle defined on $U$. 
Then the subset
$$S^{ss}(E) = \{\, s\in S \,|\, E_s \mbox{ is semistable}\,\} $$
is open in $S$.
\end{proposition}

\proof Without loss of generality, we can assume that $S$ is connected.
We proceed by induction on the relative dimension $n$ of $X/S$, starting with $n=1$.
In this case the bundle $E$ is defined on all of $X$.
By Lemma \ref{Finiteness of the set of types over S}, the set $A$ of all possible
canonical types $\HN(E_s)$ where $s\in S$ is a finite subset of $\ov{C}$.
We denote by $\tau_0 \in \ov{C}$ the type of the identity reduction of any $E_s$ to $P = G$.
This is independent of $s$ as $S$ is connected, $E_s$ is semistable if and only if
$\HN(E_s) = \tau_0$, and $\tau \ge \tau_0$ for all $\tau \in A$. 

Let $\tau \in A$ correspond to the standard parabolic $P$.  
Consider the projective morphisms $E/P \to X \to S$. 
As in Section 2 of [G-N-1],
the relative Hilbert scheme $Hilb_{(E/P)/S}$ has an open subscheme $R_{(E/P)/X/S} \subset Hilb_{(E/P)/S}$
that represents the contrafunctor which
to any $T\to S$ associated the set of all sections of $E_T/P \to X_T$.
Each such section is the same as a relative reduction of $E_T$ to $P$ over $X_T/T$.
There is
an open and closed subscheme $R^{\tau}_{(E/P)/X/S} \subset R_{(E/P)/X/S}$ that 
represents the contrafunctor which
to any $T\to S$ associated the set of relative reduction of $E_T$ to $P$ over $X_T/T$
of constant type $\tau$. It is shown in the proof of
Theorem 2.6 of [G-N-1] that $R^{\tau}_{(E/P)/X/S}$ is of finite type over $S$.

The set theoretic image $V_{\tau}\subset |S|$ of the projection $R^{\tau}_{(E/P)/X/S} \to S$
consists of all $s\in S$ such that $E_s$ admits at least one reduction of type $\tau$. 
It follows by Corollary \ref{maximality of HN type} that $\HN(E_s) \ge \tau$.
Hence as a set, 
$$S^{ss}(E) = |S| - \cup_{\tau \in A - \{ \tau_0 \}} V_{\tau},$$
which shows that $S^{ss}(E)$ is constructible. Moreover, $S^{ss}(E)$ 
is closed under generization of points
by the valuative criterion (Theorem \ref{valuative criterion}).
Hence it is open, proving the result when
the relative dimension $n$ of $X/S$ is $1$.

So now let $n\ge 2$. The proof of the inductive step is as in Proposition 7.4 of [G-N-2],
which we briefly sketch. As $S$ is noetherian,
by semicontinuity and base-change there exists an integer $m_1$ such that $\OO_{X_s}$ is $m_1$-regular 
for all $s\in S$. For any $s_0\in S$ and a smooth effective divisor $H\subset X_{s_0}$ such
that $H\in |\OO_{X_{s_0}}(m)|$ where $m\ge m_1$, there exists an open neighbourhood $s_0\in V\subset S$
and a relative effective divisor $Y\subset X_V$ over $V$ such that (i) $Y \in |\OO_{X_V/V}(m)|$
(ii) $Y\to V$ is smooth, and (iii) $Y_{s_0} = H \subset X_{s_0}$.
We can choose $H$ such that $U_{s_0}\cap H$ is a big open subscheme of $H$. By shrinking $V$ if
needed, we can assume that $U\cap Y$ is relatively big over $V$. 
Now let $s_0\in S^0(E)$. By the semistable restriction theorem 
([Gu] Theorem 12), there is an integer $m_2\ge m_1$ such that if $m\ge m_2$, then
$E_{s_0}|_{U_{s_0}\cap H}$ is semistable, defined
on a big open subscheme in $H$. 
As $Y\to V$ has relative dimension $n-1$,
by inductive hypothesis there exists a open neighbourhood $s_0\in V'\subset V$ such that
$E_s|_{Y_s}$ is semistable for all $s \in V'$. It follows by Proposition 4.4 of [G-N-2]
that $E_s$ is semistable for all $s \in V'$. Hence $V' \subset S^{ss}(E)$, which shows that 
$S^{ss}(E)$ is open. \hfill$\square$.

\bigskip

{\bf Proof of Theorem \ref{family version of main theorem}.} 
An exposition of the basic facts about 
the relative Picard scheme $Pic_{X/S}$ that we need can be found in [K].
We first treat the case where $\pi: X\to S$ admits
a global section $\sigma : S\to X$. Under this assumption, 
there exists a Poincar\'e line bundle 
on $X\times_S Pic_{X/S}$, determined up to (a non-unique) 
isomorphism by the requirement that its pullback to $Pic_{X/S}$ 
under the section 
$(\sigma, \id_{Pic_{X/S}}): Pic_{X/S} \to X\times_S Pic_{X/S}$
is a trivializable line bundle. 
The choice of a Harder-Narasimhan type $\tau \in \ov{C}$
in the closed positive Weyl chamber determines the Hilbert polynomial 
$h \in \Q[t]$ of any line bundle 
$L$ which occurs in a canonical reduction $[L, f: L\to E_s(V_P)]$
of type $\tau$ of the principal $G$-bundle $E_s = E|_{X_s}$ on a fiber $X_s$. 
Let $J = Pic_{X/S}^h \subset Pic_{X/S}$ be the open and closed subscheme where
the Hilbert polynomial of the line bundle is $h$. Let $\LL$
denote the restriction of the Poincar\'e line bundle 
to $X\times_S J\subset X\times_S Pic_{X/S}$. 
Let $\F$ be the coherent $\OO$-module on $X\times_S J$
defined by
$$\F = \un{Hom}(\LL, E_J(V_P))$$
and let $\QQ$ denote the coherent $\OO_J$-module
which is the Grothendieck Q-sheaf for $\F/X_J/J$. Let 
$Y = \PP(\QQ)$ be the corresponding projective scheme over $J$,
which represents the functor $\Psi = \Psi_{\F/X_J/J}$ by 
Lemma \ref{projective representation}. 
Over $Y$, we have a universal element 
$(M,i) \in \Psi(Y)$ as given by Lemma \ref{projective representation},
where $M = \OO_Y(-1)$ and
$i : M \to {\pi_Y}_*\F_Y$ is obtained by applying $\un{Hom}(-, \OO_Y)$
to the tautological quotient line bundle on $Y$.
By the left-adjointness of
$\pi_Y^*$ to ${\pi_Y}_*$, we have an $\OO_{X_Y}$-linear homomorphism
$f : \pi_Y^*M \to \F_Y$, which by the left-adjointness
of $\otimes$ to $\un{Hom}$ is the same as an $\OO_{X_Y}$-linear homomorphism
$$ f : \LL \otimes_{\OO_{X_Y}} \pi_Y^*M  \to E_Y(V_P).$$
Note that for any $y\in Y$, the resulting $\OO_{X_y}$-linear homomorphism
$f_y: \LL_y \otimes_{\kappa(y)} M_y \to E_y(V_P)$ is nonzero by definition of $\Psi$.

Let $W \subset X_Y$ be the open subscheme of $X_Y$ consisting of all $x\in X_y$ where 
$f_y : \LL_y \otimes_{\kappa(y)} M_y \to  E_Y(V_P)$ is injective on the fibers at $x$,
and let $X_Y - W \to Y$ be the
projection, where $X_Y - W$ is given the reduced closed subscheme structure.
By semicontinuity for dimension of fibers of a projective morphism applied
to $X_Y - W \to Y$, the subset $Y_1$ of $Y$ consisting 
of points $y$ where the fiber dimension is $\le \dim(X_y) - 2$ is open in $Y$.
We regard $Y_1$ as an open subscheme of $Y$. 
This is the open subscheme of $Y$ over which
$f : \LL\otimes_{\OO_{X_Y}} \pi_Y^*M \to  E_Y(V_P)$ is fiber wise injective 
in a relatively big open subscheme of $X_{Y_1}$ over $Y_1$.

By Remark 7.5 of [G-N-2], the set $C$ of points of $X_Y$ over
which $f : \LL\otimes \pi_Y^*M \to  E_Y(V_P)$ factors via the cone 
$\wh{E/P}\subset E_Y(V_P)$ over $E/P \subset 
{\bf P}(E_Y(V_P))$ has a natural structure of a 
closed subscheme, with an appropriate universal property. 
Consider the flattening stratification of $C \to Y$, and let $Y_2\subset Y$
be the stratum indexed by the Hilbert polynomial of the fiber of $X/S$.
This is the top stratum so $Y_2$ is a closed subscheme of $Y$.
It has the universal property that under a base change $T\to Y$ the pullback
of $f$ factors via the pullback of the cone $\wh{E/P}$ if and only if
the morphism $T\to Y$ factors via $Y_2\subset Y$.

Let $Y_3 = Y_1\cap Y_2 \subset Y$, which is the locally closed
subscheme of $Y$ over which $f : \LL \otimes_{\OO_{X_Y}}  \pi_Y^*M  \to  E_Y(V_P)$ defines
a rational reduction of structure group to $P\subset G$. 
Let $Y_4\subset Y_3$ be the open and closed subscheme of $Y_3$ where the 
type of the reduction defined by $f$ is $\tau$. Finally, by Proposition \ref{openness of semistability},
we have an open subscheme 
$Y_5 \subset Y_4$ where the extension under
the Levi quotient $P \to P/R_u(P)$ is semistable. 
It is immediate from its construction that the $S$-scheme $Z = Y_5$ 
represents the functor $T\mapsto \Phi_{E/X/S}^{\tau}(T)$ which is the set of all
relative canonical reductions of type $\tau$ of the pullback $E_T/X_T/T$.

Now we come to a general case, where $X$ may not necessarily admit a 
global section over $S$. As $X\to S$ is by assumption smooth,
there exists a finite-type surjective separated \'etale morphism $p : S'\to S$ 
such that the base change $X'= X_{S'}$ admits a global section
$S' \to X'$. Let $E' = E_{S'}$. 
Hence by the above special case, there exists a scheme 
$Z'\to S'$ which represents the functor 
$\Phi_{E'/X'/S'}^{\tau} : (Schemes/S')^{opp} \to Sets$, which is just
the functor $p^*\Phi_{E/X/S}^{\tau} = (\Phi_{E/X/S}^{\tau})_{/S'}$.
The functor $\Phi_{E/X/S}^{\tau}$ is an fppf sheaf by Lemma \ref{sheaf property}.
It therefore follows by
Lemma \ref{descent of representating space} below, that there exists
an algebraic space $Z \to S$ which represents $\Phi_{E/X/S}^{\tau}$,
and moreover, we have an isomorphism $Z\times_S S' = Z'$.
As $Z'\to S'$ is of finite type and separated, 
it follows by fppf descent that $Z \to S$ is also so.

Next, we note that the uniqueness of a canonical reduction in the absolute case
(Theorem \ref{existence and uniqueness})
implies that $Z'\to S'$ is 
injective at the level of underlying sets. By the above
reasoning, this implies that $Z\to S$ too is injective
at the level of underlying sets. As this holds for all base changes,
it follows that $Z\to S$ universally injective (radicial).
As $Z\to S$ is of finite type and injective, it follows that
$Z\to S$ is quasi-finite.

Given the above properties of $Z\to S$, 
the Proposition \ref{SP representability by scheme} below
implies that $Z$ is a scheme. This is the desired scheme $S^{\tau}(E)$ by
its construction. It is of finite type over $S^{\tau}(E)$ by construction.

If $x\in S^{\tau}(E)$ maps to $y \in S$, consider the bundle $E_y$ on
$X_y = X \times_S \spec k(y)$, which is equipped with a unique canonical reduction
and it has type $\tau$. Hence by the universal property of $S^{\tau}(E)\to S$,
the morphism $y: \spec k(y) \to S$ had a unique lift $y': \spec k(y) \to S^{\tau}(E)$.
As $x$ is the only point in the fiber of $S^{\tau}(E)\to S$ over $y$, it follows
that the image of $y': \spec k(y) \to S^{\tau}(E)$ is $x$, and $k(y) \to k(x)$
is an isomorphism. This proves the assertion of the
Theorem \ref{family version of main theorem} that the morphism
$S^{\tau}(E)\to S$ induces an isomorphism on all residue fields of points
of $S^{\tau}(E)$. 

Note that for any $\tau\in \ov{C}$,
the set-theoretic image of $S^{\tau}(E)\to S$ is the subset
$$|S|_{\tau} = \{ s\in S\,|\, \HN(E_s) = \tau\}$$
which is therefore a constructible subset of $S$ by
Chevalley's constructibility theorem.

By Proposition 
\ref{Finiteness of the set of types over S}, the set
$$A = \{ \tau \in \ov{C}\,|\, \tau = \HN(E_s)\mbox{ for some }s\in S\}$$
is finite. 
Taking a finite union, this shows that for any $\tau \in \ov{C}$, the subset 
$$|S|_{> \tau} = \{ s\in  S\,|\, \HN(E_s) > \tau\} \subset |S|$$
is a constructible subset of $|S|$, where $|S|$ denotes the underlying subset of
$S$.  By an application of
Theorem \ref{valuative criterion} and Lemma \ref{exercise}, 
this constructible subset is closed in $S$, 
hence its complement 
$$S^{\ngtr \tau}(E) = S - |S|_{> \tau}\,  = \{ s\in S \,|\, \HN(E_s) \ngtr \tau \}
\subset S$$
is open in $S$. We give it the structure of an 
open subscheme of $S$. The morphism $S^{\tau}(E)\to S$ factors via
the inclusion $S^{\ngtr \tau}(E)\hra S$, inducing a morphism 
$S^{\tau}(E)\to S^{\ngtr \tau}(E)$. This is a morphism of finite type between
noetherian schemes, and by Theorem \ref{valuative criterion} it satisfies the 
valuative criterion for properness w.r.t. discrete valuation rings. Hence
the induced morphism $S^{\tau}(E)\to S^{\ngtr \tau}(E)$ is proper.

In case the Behrend conjecture is satisfied, $S^{\tau}(E)\to S^{\ngtr \tau}(E)$ is
unramified by the infinitesimal uniqueness of a canonical reduction. An unramified
radicial finite morphism is a closed embedding, which proves the last assertion
in the statement of the theorem. 

This completes the proof of the Theorem \ref{family version of main theorem},
and also of the Theorem \ref{semicontinuity of type} on semicontinuity of
the HN-type, stated below.
\hfill$\square$

\begin{theorem}\label{semicontinuity of type}{\bf Semicontinuity of
    canonical type}
  Let $S$ be a noetherian scheme, and let
  $\pi: X\to S$ be a smooth projective morphism with geometrically
  irreducible fibers, together with a relatively very ample line bundle $\OO_{X/S}(1)$.
  Then for any principal $G$-bundle $E$ on $X$, the following is true for
  any $\tau\in \ov{C}$.

(1) The subset
$|S|_{\tau}(E)    =  \{ s\in S \,|\, \HN(E_s) = \tau \}$ is a locally closed subset of
$|S|$.

(2) The closure of $|S|_{\tau}(E)$ is contained in the union
$\cup_{\sigma \ge \tau}\, |S|_{\sigma}(E)$, which is a closed subset of $|S|$.

(3) The subset $|S|_{> \tau}(E) = \cup_{\sigma > \tau}\, |S|_{\sigma}(E)$
is a closed subset of $|S|$, equivalently, 
$S^{\ngtr \tau}(E) =  \{ s\in S \,|\, \HN(E_s) \ngtr \tau \} \subset S$
is open in $S$. The subset $|S|_{\tau}(E)$ is closed in $S^{\ngtr \tau}(E)$.
\hfill$\square$
\end{theorem}  

Note that for $\tau = \mu_{(G,\id)}$ (the type of the identity reduction to $P=G$),
we have $|S|_{\tau}(E) = S^{ss}(E)$ and in that case 
the above assertions (1), (2) and (3) follow from the Proposition \ref{openness of semistability}
which more generally applies to bundles defined on big opens.
The Proposition \ref{openness of semistability},  
applied to Levi reductions (which may only be defined on big opens) was used in the proof
of the above theorem also for a general $\tau$.

\begin{lemma}\label{descent of representating space}
Let $S$ be a scheme and $\Phi: (Schemes/S)^{op} \to Sets$ a functor which
is an fppf sheaf of sets on $S$.
Suppose that there exists an fppf morphism $p : S' \to S$ such that  
the functor $p^*\Phi = \Phi_{/S''} : (Schemes/S')^{op} \to Sets$ is representable
by an algebraic space $Z'/S'$. Then the functor
$\Phi$ is representable by an algebraic space $Z/S$.
\end{lemma}

\proof 
Let $S'' = S'\times_S S'$, and let 
$p_1,p_2 : S'' \stackrel{\to}{\scriptstyle\to} S'$ be the 
two projections. We write $\pi = p\circ p_1 = p\circ p_2$,
so that we have natural isomorphisms
$$p_1^*(\Phi_{/S''}) =  \pi^*\Phi = p_2^*(\Phi_{/S''}).$$ 
Hence we have an isomorphism of the representing $S''$-spaces
$p_1^*Z' \to p_2^*Z'$ which satisfies the cocycle condition.
As fppf descent is effective for algebraic spaces
(see Lemma 78.11.3 of [Stacks Project]), it follows that 
$Z'$ descends to define an algebraic space $Z$ over $S$,
whose pull back under $S'\to S$ is $Z'$. 
\hfill$\square$

Actually, we only need the case of the above lemma where $S'\to S$ is \'etale,
and so effectiveness of \'etale descent for algebraic spaces (see
Corollary 1.6.4 of [La-MB]) is enough for us.

\begin{proposition}\label{SP representability by scheme} 
{\rm ([Stacks Project] Tag 03XX, Proposition 55.47.2.)} 
Let $S$ be a scheme. Let $f:X \to T$ be a morphism of 
algebraic spaces over $S$. Assume that 
$T$ is representable by a scheme, and 
$f$ is separated and locally quasi-finite. 
Then $X$ is representable by a scheme. 
\end{proposition}

{\bf Proof of Theorem \ref{stack version of main theorem}.} 
We begin by recalling that the  
stack $Bun_{X/S}(G)$ of $G$-bundles on fibers of $X/S$ is algebraic. 
To see this, choose a closed embedding $G \hra GL_{n,k}$
as group schemes over $k$, and consider
the induced $1$-morphism of stacks $Bun_{X/S}(G)\to Bun_{X/S}(GL_{n,k})$. 
The stack $Bun_{X/S}(GL_{n,k})$ is just the 
stack of rank $n$ vector bundles on fibers of $X/S$, so it is an 
algebraic stack (see [La-MB] Theorem 4.6.2.1). Given any $GL_{n,k}$-bundle $E$ on $X$,
the reductions of its structure group to $G$ are the sections of 
$E/G \to X$, so they are parameterized by a suitable open subscheme of  
the Hilbert scheme $Hilb_{(E/G)/S}$ (see for example [N-2] section 5.6.2
for an exposition). 
This shows the $1$-morphism $Bun_{X/S}(G)\to Bun_{X/S}(GL_{n,k})$ is schematic, 
which implies that the stack $Bun_{X/S}(G)$ is algebraic.

Next, given $\tau \in \ov{C}$, consider the $1$-morphism from the stack 
$Bun_{X/S}^{\tau}(G)$  of Corollary 1.2 to the stack $Bun_{X/S}(G)$. 
The Theorem \ref{family version of main theorem} 
shows that this $1$-morphism is schematic and has the desired properties.
\hfill$\square$

{\bf Acknowledgement} S. Gurjar thanks the Tata Institute of 
Fundamental Research, and N. Nitsure thanks
the Indian Institute of Technology - Bombay,
for their hospitality and support
during a part of the preparation of this article.

\bigskip

\bigskip

{\footnotesize

\parskip=2pt
{\large \bf References}

\bigskip

[Be] Behrend, K. : Semi-stability of reductive group schemes over curves. 
Math. Ann. 301 (1995), 281-305.

[Bi-Ho] Biswas, I. and Holla, Y. I. : Harder-Narasimhan reduction of a principal 
bundle. Nagoya Math. J. 174 (2004), 201-223

[EGA] Grothendieck, A. and Dieudonn\'e, J. : 
{\it  \'El\'ements de g\'eom\'etrie alg\'ebrique}, 
Publ. Math. IHES., vols. 4, 8, 11, 17, 20, 24, 28, 32 (1960-1967). 

[Gu] Gurjar, S. : 
Restriction theorems for principal bundles in arbitrary characteristic. 
J. Algebra 426 (2015), 79-91. 

[G-N-1] Gurjar, S. and Nitsure, N. : Schematic Harder-Narasimhan 
stratification for families of principal bundles and $\Lambda$-modules. 
Proc. Indian Acad. Sci. (Math. Sci.) 124 (2014), 315-332. 


[G-N-2] Gurjar, S. and Nitsure, N. : 
Schematic Harder-Narasimhan stratification for families of 
principal bundles in higher dimensions. Math. Zeit. 289 (2018), 1121-1142.

[He] Heinloth, J. : Bounds for Behrend's conjecture on the canonical 
reduction. Int. Math. Res. Not. IMRN Vol. 2008, rnn045.

[La-MB] Laumon, G. and Moret-Bailly, L. : {\it Champs alg\'ebriques},
Springer (2000).

[K] Kleiman, S. : The Picard scheme. Part 5 
of {\it Fundamental Algebraic Geometry -- Grothendieck's FGA 
Explained}, Fantechi et al, Math. Surveys and Monographs Vol.
123, American Math. Soc. (2005).

[M-R] Mehta, V. B. and Ramanathan, A. : Semistable sheaves on 
projective varieties and their restriction to curves.
Math. Ann. 258 (1981/82), 213-224.

[M] Milne, J. S. : {\it Algebraic Groups}, Cambridge University Press (2017).

[N-1] Nitsure, N. : Representability of Hom implies flatness. 
Proc. Indian Acad. Sci. Math. Sci. {\bf 114} (2004)

[N-2] Nitsure, N. : Construction of Hilbert and Quot schemes. Part 2 
of {\it Fundamental Algebraic Geometry -- Grothendieck's FGA
Explained},  Fantechi et al, Math. Surveys and Monographs Vol.
123, American Math. Soc. (2005).

[N-3] Nitsure, N. : 
Schematic Harder-Narasimhan stratification. 
Internat. J. Math. 22 (2011), 1365-1373. 

[Rag] Raghunathan, M. S. : A note on quotients of real algebraic groups by
arithmetic subgroups. Invent. Math. 4 (1967/1968), 318-335. 

[Ram] Ramanathan, A. : Moduli for principal bundles over algebraic curves, I and II.
Proc. Indian Acad. Sci. (Math. Sci.). vol 106 (1996). Part I: pages 301-328.
Part II: pages 421-449.

[Stacks Project] : http://stacks.math.columbia.edu/

\bigskip

\bigskip



Sudarshan Gurjar \hfill Nitin Nitsure\\
Department of Mathematics  \hfill Professor of Mathematics (retired)\\
Indian Institute for Technology - Bombay 
\hfill Tata Institute of Fundamental Research\\
Mumbai 400 076 \hfill Mumbai 400 005\\ 
India \hfill India\\
{\tt srgurjar1984@gmail.com} \hfill {\tt nitsure@gmail.com}

}

\end{document}